\documentclass[12pt]{iopart}
\usepackage[latin1]{inputenc}
\usepackage{iopams}
\usepackage{harvard}
\usepackage{graphicx}
\usepackage{psfrag}

\newtheorem{theorem}{Theorem}[section]
\newtheorem{lemma}[theorem]{Lemma}

\newtheorem{corollary}[theorem]{Corollary}
\newtheorem{remark}[theorem]{Remark}

\newtheorem{definition}[theorem]{Definition}
\newtheorem{conjecture}[theorem]{Conjecture}

\newcommand{\field}[1]{\mathbb{#1}} 
\def\C{\mathbb{C}}
\def\R{\mathbb{R}}

\def\fC{\field{C}}
\def\fR{\field{R}}
\def\fZ{\field{Z}}

\def\cA{\mathcal{A}}
\def\cB{\mathcal{B}}
\def\cC{\mathcal{C}}
\def\cD{\mathcal{D}}
\def\cR{\mathcal{R}}
\def\cP{\mathcal{P}}
\def\cO{\mathcal{O}}

\def\bF{{\bf F_*}}
\def\bs{{\bf s^*}}
\def\loc{{\rm loc}}

\date{2011-06-02}

\begin{document}
\title{No elliptic islands for the universal area-preserving map}
\author{Tomas Johnson}
\address{Department of Mathematics, Cornell University, Ithaca, NY 14853, USA}
\eads{\mailto{tomas.johnson@cornell.edu}}
\begin{abstract}

A renormalization approach has been used in \cite{EKW1} and \cite{EKW2} to prove the existence of a \textit{universal area-preserving map}, a map with hyperbolic orbits of all binary periods. The existence of a horseshoe, with positive Hausdorff dimension, in its domain was demonstrated in \cite{GJ1}. In this paper the coexistence problem is studied, and a computer-aided proof is given that no elliptic islands with period less than $20$ exist in the domain. It is also shown that less than $1.5\%$ of the measure of the domain consists of elliptic islands. This is proven by showing that the measure of initial conditions that escape to infinity is at least $98.5\%$ of the measure of the domain, and we conjecture that the escaping set has full measure. This is highly unexpected, since generically it is believed that for conservative systems hyperbolicity and ellipticity coexist.
\end{abstract}
\submitto{Nonlinearity}
\maketitle
\setcounter{page}{1}


\section{Introduction}
The phenomena of \textit{coexistence}, simultaneous existence of hyperbolicity, often formulated as the existence of positive Lyapunov exponents, and stability, in the form of elliptic islands, is believed to be generic in Hamiltonian dynamics \cite{B08,G10,HH64,L04,MM74,P82,S91}. Very few examples of systems where coexistence has been rigorously proved exist; the first such example is given in \cite{P82}. At the same time very few, if any, systems exist that numerically do not exhibit coexistence, as soon as hyperbolicity is present.

In this paper we demonstrate that the \textit{universal} area-preserving map associated with period doubling does not have any elliptic islands of low periods, and that the area of the remaining ones is small. This is a non-trivial result, since for the universal map, it is only known that one family of elliptic islands of period $2^n$ have bifurcated into hyperbolic periodic orbits. Nothing is known apriori about the periodic orbits with period different from $2^n$, and the orbits with period $2^n$ different from those involved in the bifurcation sequence. The computations give numerical evidence towards the nonexistence of elliptic islands for the universal area-preserving map. Together with our resent paper \cite{GJ1}, in which we prove that the same map has a hyperbolic horseshoe with positive Hausdorff dimension, this gives a candidate for an example of a conservative system that does not exhibit coexistence. In fact, these results are likely to hold for all infinitely renormalizable maps in a neighborhood of the universal one. For area-preserving H\'enon maps, this contrasts to parameters close to the beginning of the period doubling sequence, where the hyperbolic sets are accumulated by elliptic islands as proved in \cite{Duarte2}.

Our proof is computer-aided and what we actually prove is that there are no elliptic islands of period less than $20$, and that the area enclosed by the elliptic islands is less than $0.049$, which is less than $1.5\%$ of the domain. The estimate of the bound on the area enclosed by the invariant set scales down like the area of $n$-th generation covers of a Cantor set, indicating that the area actually enclosed by elliptic islands is zero. It should be noted, however, that the computations are done using the bound on the fixed point from \cite{EKW2}, and most maps in this neighborhood have elliptic islands. Since what we compute are successive approximations of the invariant set, it appears that the measure of the invariant set is zero, i.e., that the set of initial conditions that escape to infinity has full measure. This is in agreement with the approximations of the measure of the invariant set of the conservative H\'enon family computed in \cite{M82}. 

The universal area-preserving map is a map with hyperbolic periodic points of all binary periods $2^n$. Its existence, as many other universality phenomena is best described in the setting of renormalization. Universality will stand for the independence of the quantifiers of the geometry of orbits and bifurcation cascades in families of maps of the choice of a particular family. To prove universality one usually introduces a {\it renormalization} operator on a functional space, and demonstrates that this operator has a hyperbolic fixed point. Universality was first discovered for the case of period doubling universality in unimodal maps - the Feigenbaum-Coullet-Tresser universality \cite{Fei1,Fei2,TC}. Proofs appear in e.g. \cite{Eps1,Eps2,Lyu}.

It has been established that the universal behavior in dissipative and conservative higher dimensional systems is fundamentally different. The case of the dissipative systems is often reducible to the one-dimensional Feigenbaum-Coullet-Tresser universality \cite{CEK1,dCLM,LM,HLM10}. The case of area-preserving maps seems to be very different, and at present there is no deep understanding of universality in conservative systems, other than in the case of universality for systems near integrability \cite{Koch1,Koch2,Gai1,Kocic,KLDM}.

An infinite period-doubling cascade in families of area-preserving maps was observed by several authors in the early 80's \cite{DP,Hel,BCGG,Bou,CEK2,M93}. This universality can be explained rigorously if one shows that the {\it renormalization} operator
\begin{equation}\label{Ren}
R[F]=\Lambda^{-1}_F \circ F \circ F \circ \Lambda_F,
\end{equation}
where $\Lambda_F$ is some $F$-dependent coordinate transformation, has a fixed point, and the derivative of this operator is hyperbolic at this fixed point.

It has been argued in \cite{CEK2}  that $\Lambda_F$ is a diagonal linear transformation. Furthermore, such $\Lambda_F$ has been used in \cite{EKW1} and \cite{EKW2} in a computer assisted proof of existence of a reversible renormalization fixed point $F_*$ and hyperbolicity of the operator $R$. We began the study of the dynamics of the renormalization fixed point $F_*$ in \cite{GJ1,GJ2,GJM}. Additional properties of the universality phenomena in the area-preserving case are studied in \cite{GK1}, where it is demonstrated that the fixed point $F_*$ is very close, in some appropriate sense, to an area-preserving H\'enon-like map
\begin{equation*}
H^*(x,u)=(\phi(x)-u,x-\phi(\phi(x)-u )),
\end{equation*}
where $\phi$ solves the following one-dimensional problem of non-Feigenbaum type:
\begin{equation*}\label{0_equation}
\phi(y)={2 \over \lambda} \phi(\phi(\lambda y)) -y.
\end{equation*}

In \cite{GJ1} we constructed a hyperbolic set, using the method of covering relations \cite{ZG04,Z09} in rigorous computations. To prove that the constructed horseshoe has positive Hausdorff dimension we used the Duarte Distortion Theorem \cite{Duarte1,GoKa} which enables one to use the distortion of a Cantor set to find bounds on the dimension. 

The structure of this paper is as follows: in Section 2 we recall the necessary notation and concepts about area-preserving maps, their renormalization theory and the description of the domain of $F_*$, in Section 3 we summarize the results from \cite{GJ1} about the hyperbolicity of $F_*$, in Section 4 we state our theorems and conjectures, in Section 5 we describe our computer-aided proof, and finally in Section 6 we describe the results of our computations.  


\section{Renormalization of reversible area-preserving maps} 
An \textit{area-preserving map} will mean an exact symplectic diffeomorphism of a subset of ${\fR}^2$ onto its image.

Recall that an area-preserving map can be uniquely specified by its generating function $S$ (under the assumption that $(x,y)$ is a coordinate system):
\begin{equation}\label{gen_func}
\left( x \atop -S_1(x,y) \right) {{ \mbox{{\small \it  F}} \atop \mapsto} \atop \phantom{\mbox{\tiny .}}}  \left( y \atop S_2(x,y) \right), \quad S_i \equiv \partial_i S.
\end{equation}

Furthermore, we will assume that $F$ is reversible, that is 
\begin{equation}\label{reversible}
T \circ F \circ T=F^{-1}, \quad {\rm where} \quad T(x,u)=(x,-u).
\end{equation}

For such maps it follows from $(\ref{gen_func})$ that 
\begin{equation}\label{littles}
S_1(y,x)=S_2(x,y) \equiv s(x,y), 
\end{equation}

and
\begin{equation}\label{sdef}
\left({x  \atop  -s(y,x)} \right)  {{ \mbox{{\small \it  F}} \atop \mapsto} \atop \phantom{\mbox{\tiny .}}} \left({y \atop s(x,y) }\right).
\end{equation}

It is this little $s$ that will be referred to below as \textit{the generating function}. We use the same normalization as in \cite{EKW2}, i.e., $s(1,0)=0$ and $\partial_2s(1,0)=0.2$. It follows from (\ref{littles}) that $s_1$ is symmetric. If the equation $-s(y,x)=u$ has a unique differentiable solution $y=y(x,u)$, then the derivative of such a map $F$ is given by the following formula:

\begin{equation}\label{Fder}
\hspace{-1.5cm}DF(x,u)=\left[ 
\begin{array}{c c}
-{s_2(y(x,u),x) \over s_1(y(x,u),x)} &  -{1 \over s_1(y(x,u),x)} \\
s_1(x,y(x,u))-s_2(x,y(x,u)) {s_2(y(x,u),x) \over s_1(y(x,u),x)}  & -{s_2(x,y(x,u)) \over s_1(y(x,u),x)} 
\end{array}
\right]. 
\end{equation}

The following definition of the renormalization operator acting on generating functions is due to the authors of \cite{EKW1} and \cite{EKW2}:

\medskip

\begin{definition}
\phantom{a}

\begin{eqnarray}\label{ren_eq}
\nonumber \\ {\cR}_{EKW}[s](x,y)=\mu^{-1} s(z(x,y),\lambda y),
\end{eqnarray}
where
\begin{eqnarray*}
\label{midpoint_eq} 0&=&s(\lambda x, z(x,y))+s(\lambda y, z(x,y)), \\
0&=&s(\lambda,1)+s(0,1) \quad {\rm and} \quad \mu=\partial_1 z (1,0).
\end{eqnarray*}
\end{definition}

\medskip

\begin{definition}
The Banach space of functions  $s(x,y)=\sum_{i,j=0}^{\infty}c_{i j} (x-0.5)^i (y-0.5)^j$, analytic on a bi-disk
$$ |x-0.5|<\rho, |y-0.5|<\rho,$$
for which the norm
$$\|s\|_\rho=\sum_{i,j=0}^{\infty}|c_{i j}|\rho^{i+j}$$
is finite, will be referred to as $\cA(\rho)$.
$\cA_s(\rho)$ will denote its symmetric subspace $\{s\in\cA(\rho) : s_1(x,y)=s_1(y,x)\}$.
\end{definition}

\medskip

As we have already mentioned, the following has been proved with the help of a computer in \cite{EKW1} and \cite{EKW2}:
\begin{theorem}\label{EKWTheorem}
There exist a polynomial $s_{\rm a} \in \cA_s(\rho)$ and  a ball $\cB_r(s_{\rm a}) \subset \cA_s(\rho)$, $r=6.0 \times 10^{-7}$, $\rho=1.6$, such that the operator ${\cR}_{EKW}$ is well-defined, analytic and compact on $\cB_r$. 

Furthermore, its derivative $D {\cR}_{EKW} \arrowvert_{\cB_r }$ has exactly two eigenvalues $\delta_1$ and $\delta_2$ of modulus larger than $1$, while 
$${\rm spec}(D {\cR}_{EKW} \arrowvert_{\cB_r }) \setminus \{\delta_1,\delta_2 \} \subset \{z \in \C: |z| \le \nu < 1\}.$$ 

Finally, there is an $s^* \in \cB_r$ such that
$$\cR_{EKW}[s^*]=s^*.$$
The scalings $\lambda_*$ and $\mu_*$ corresponding to the fixed point $s^*$ satisfy
\begin{eqnarray*}
\label{lambda} \lambda_* \in [-0.24887681,-0.24887376], \\
\label{mu} \mu_* \in [0.061107811, 0.061112465].
\end{eqnarray*}
 \end{theorem}

\medskip
 \begin{remark}
The radius of the contracting part of the spectrum ${\rm spec}(D {\cR}_{EKW}(s_*))  \setminus \{\delta_1,\delta_2 \}$ has been estimated in \cite{EKW2} to be $\nu=0.8$. In \cite{GJM} it is proved that $\nu = |\lambda_*|$.
 \end{remark}
It follows from the above theorem that there exist a codimension $2$ local stable manifold $W^s_\loc(s^*)\subset \cB_r$.

The bound on the fixed point generating function $s^*$ will be called $\bs$:
\begin{equation*}\label{bs}
\bs \equiv \left\{ s \in \cA_s(\rho) : \|s-s_{\rm a}\|_{\rho} \le  r=6.0 \times 10^{-7}  \right\},
\end{equation*}
while the bound on the renormalization fixed point $F_*$ will be referred to as  $\bF$:
\begin{equation*}\label{bF}
\bF \equiv \left\{ F: (x,-s(y,x)) \mapsto (y,s(x,y)) : s \in \bs \right\},
\end{equation*}
where $s_{\rm a}$ is as in Theorem $\ref{EKWTheorem}$.

We will now summarize the description of the domain of the universal map, $F_*$, from \cite{GJ1}. $F_*$ is defined implicitly by the generating function $s_*$, whose domain is given as in \cite{EKW2}:
\begin{equation*}\label{sDom}
\cD_s=\{(x,y) \in \fC^2 \,: \, |x-0.5|<1.6, |y-0.5|<1.6\}.
\end{equation*}
To compute the domain of $F\in\bF$, we note that its second argument is equal to $-s(y,x)$, for some $s\in\bs$, see $(\ref{sdef})$. 

Thus, the domain of $F$, $\cD$, is given by:
\begin{equation*}\label{fDom}
\!\!\!\!\!\!\!\!\!\!\!\!\!\!\!\!\!\!\!\!\!\!\!\!\!\!\!\!\!\!\!\!\! \cD = \{(x,u) \in \fC^2 \,:\, u=-s(y,x), \, s \in  \bs(y,x), \, |x-0.5|<1.6,\, |y-0.5|<1.6 \}.
\end{equation*}

\begin{figure}[h]
\begin{center}
\includegraphics[width=0.8\textwidth]{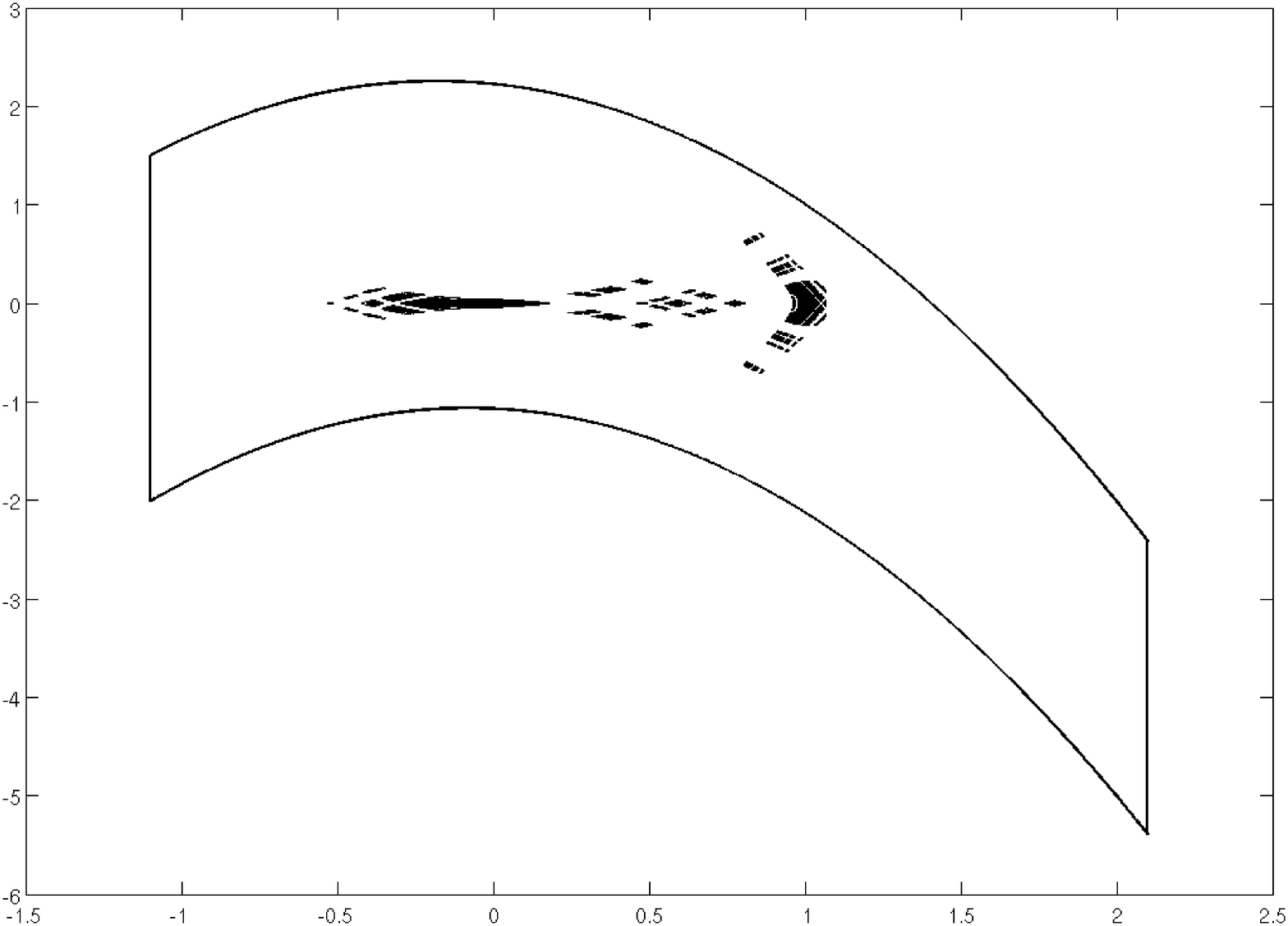}
\caption{The real slice of the domain of $F_*$ together with its invariant subset.}\label{DomPic}
\end{center}
 \end{figure}

We denote by
\begin{equation*}\label{fDomRe}
\tilde{\cD} =\{(x,u) \in \cD \,: \Im{x}=\Im{u}=0\},
\end{equation*}
the real slice of $\cD$. In \cite{GJ1} we constructed a nonempty subset of $\tilde{\cD}$, which is described in the following lemma and illustrated in Figure \ref{DomPic}.

\begin{lemma}\label{lfDomRe}
There exists a non-empty open set $\bar{\cD}\subset\R^2$, 
$$\cP_x \bar{\cD} \subset \{x \in \fR: |x-0.5| <1.6\},$$
such that for every $(x,u) \in \bar{\cD}$ and every $s \in \bs(y,x)$ there exists a unique real solution of the equation $u=-s(y,x)$ that satisfies $|y-0.5|<1.6$.
\end{lemma}

\begin{figure}[h]
\begin{center}
\includegraphics[width=0.32\textwidth]{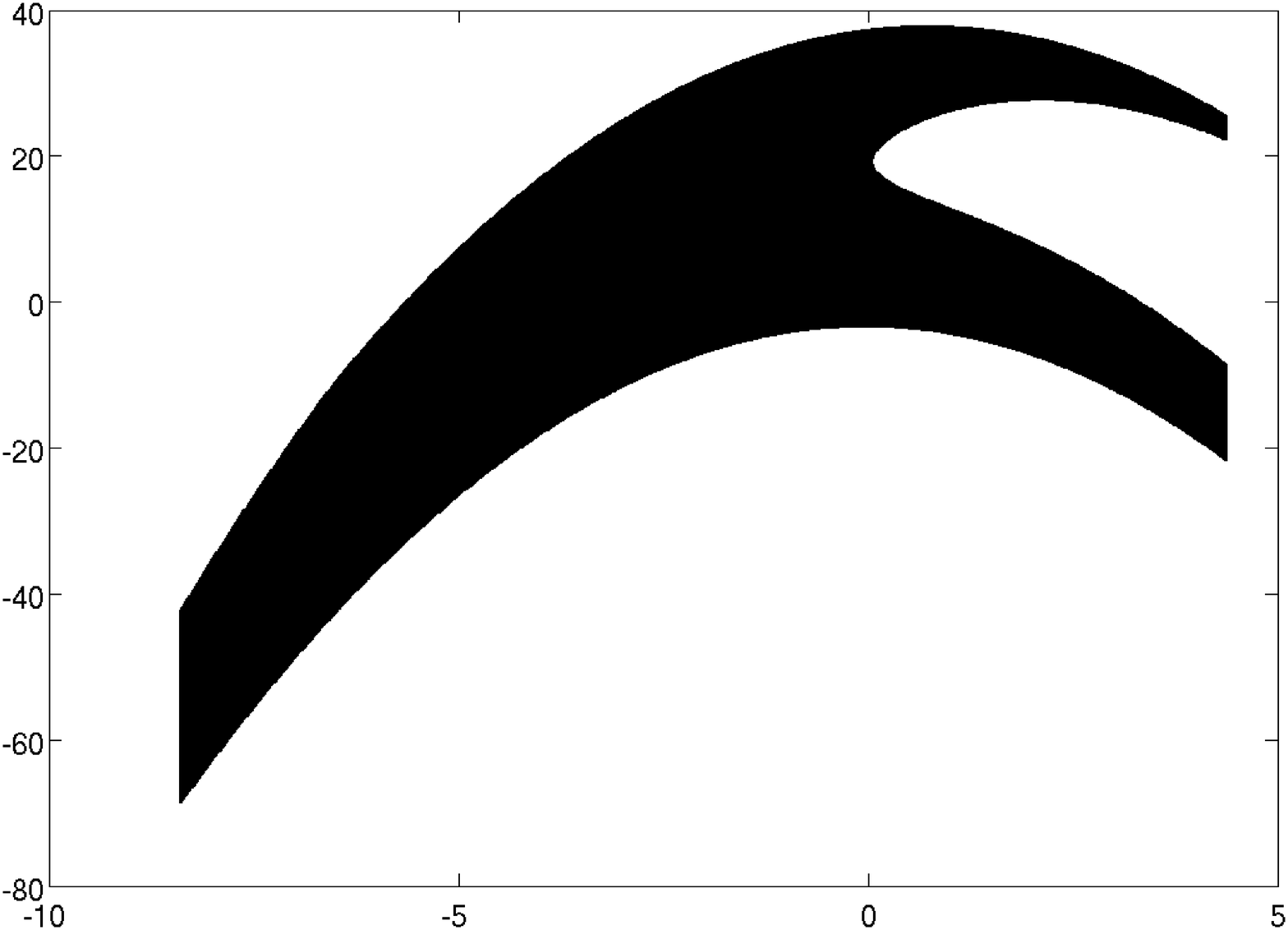}
\includegraphics[width=0.32\textwidth]{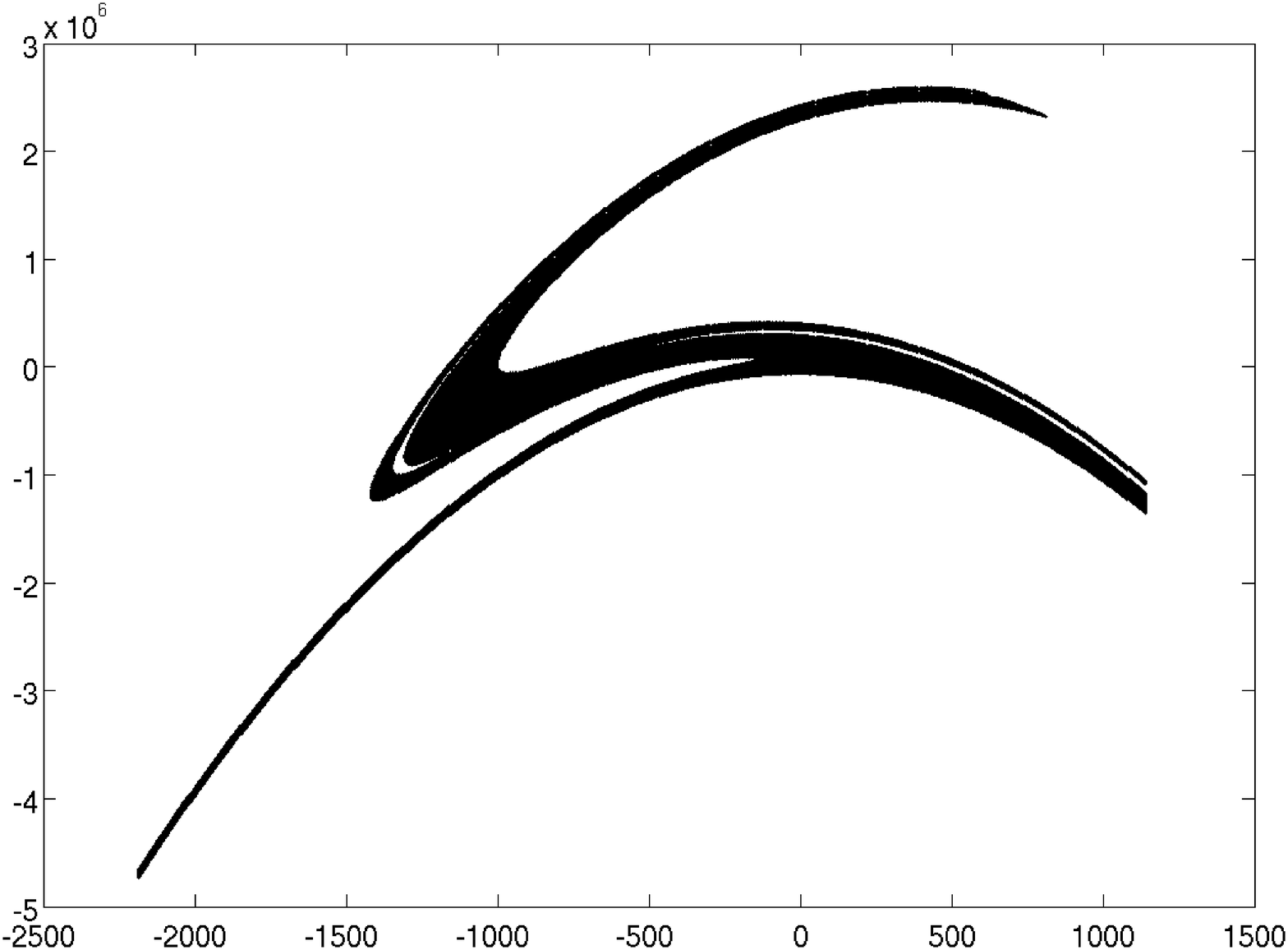}
\includegraphics[width=0.32\textwidth]{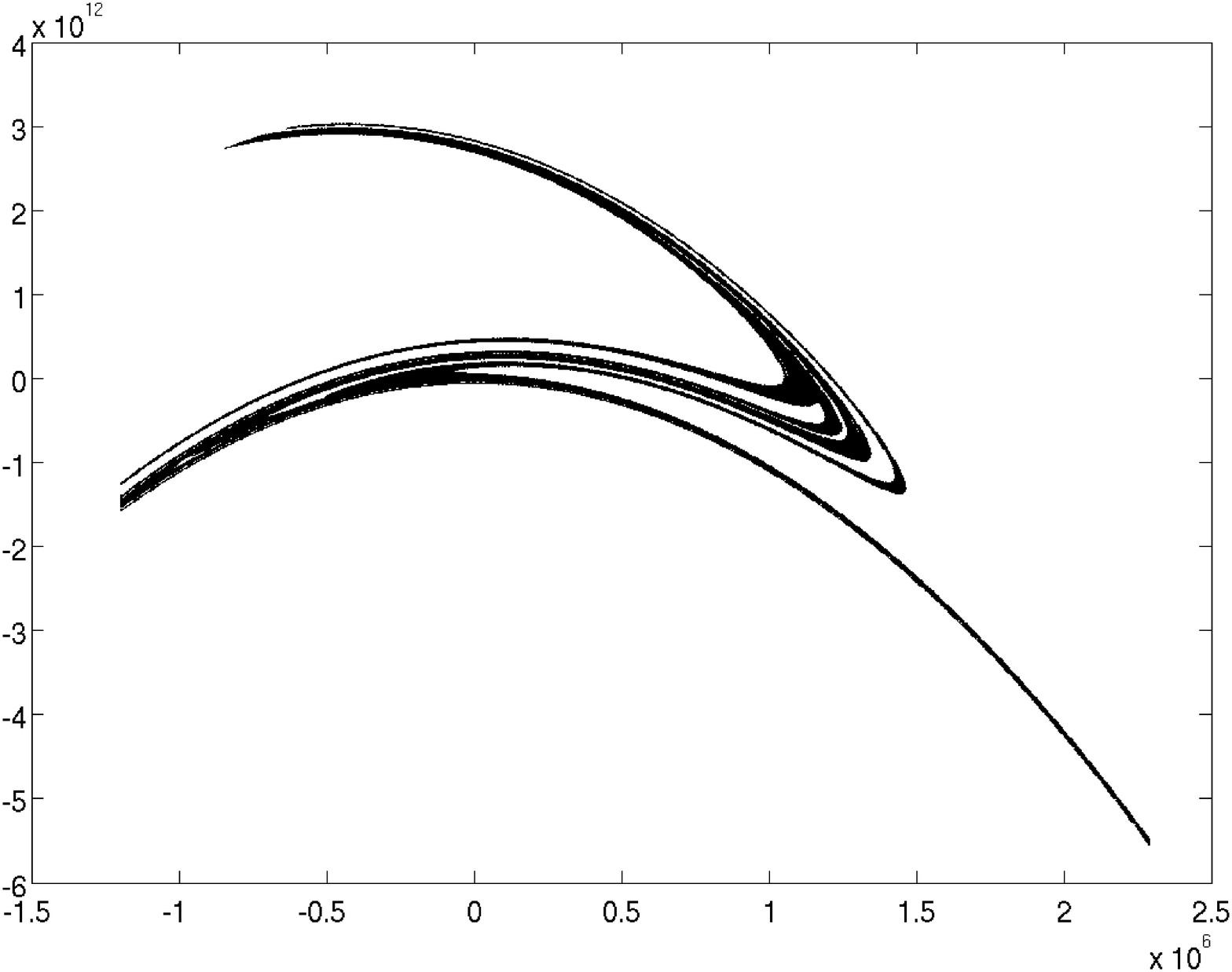}
\caption{The real slices of the domains $\cD^1, \cD^5,$ and $\cD^{10}$ of three analytic continuations of $F_*$.}\label{fD123}
\end{center}
 \end{figure}

Clearly, $\bar{\cD} \subset \tilde{\cD}$. It also follows that $\cD$ contains an open complex neighborhood of the set $\bar{\cD}$. The measure of the set $\bar\cD$ is at least $3.312$. With the invariant subset of $\cD$, we denote the set of points whose orbits stay in $\tilde \cD$. We note that by using the renormalization equation $F_*=R[F_*]$, $F_*$ has, for any positive integer $k$, an analytic continuation to the domains 
\begin{equation}
\cD^{k} := \Lambda_*^{-1} (F_*^{-1}(F_*(\cD^{k-1})\cap  \cD^{k-1})), \quad \cD^0 := \cD,
\end{equation} 
where we note that each $\cD^k$ is the rescaling of the subset of a simply connected set that is mapped into itself by $F_*$, which is an orientation preserving diffeomorphism. If $F_*(\cD^{k-1})\cap \cD^{k-1}$ is not connected, we restrict to the connected component containing the point $p_0$ described in Theorem \ref{MTHM1}. This set is simply connected, since $D^{k-1}$ and $F_*(\cD^{k-1})$ are. Thus, inductively, for each $k$, $\cD^k$ is simply connected. We put
\begin{equation}
\cD^{\infty} := \bigcup_{k=0}^\infty \cD^k.
\end{equation} 
The real slices of the domains are denoted by $\tilde \cD^k$. The three analytic continuations of $F_*$ to $\tilde \cD^1, \tilde\cD^5,$ and $\tilde\cD^{10}$ are shown in Figure \ref{fD123}. We note that for all the approximations we have computed, $F_*(\tilde \cD^{k-1})\cap \tilde \cD^{k-1}$ is connected and $\tilde \cD^{k-1}\subset\tilde \cD^k$. Similarly, $F_*^{-1}$ has analytic continuations to the sets $T\cD^k$. 

\section{Review of hyperbolic behavior}
We will now summarize our findings from \cite{GJ1} that are needed to state the result, the wording has been slightly modified to better suit the present paper:

\begin{theorem}\label{MTHM1}
The renormalization fixed point $F_*$ has the following properties:
\begin{itemize}
\item[1)] $F_*$ possesses a hyperbolic fixed point $p_0$;
\item[2)] $F_*$ possesses a point $p_\pitchfork$ which is transversally homoclinic to the fixed point $p_0$; 
\item[3)] there exists a positive integer $n$ such that for any negative integer $k$ the map $F_*^n$ has a heteroclinic orbit $\cO_k$ between the periodic points  $\Lambda_*^{k}(p_0)$ and $\Lambda_*^{k+1}(p_0)$, and for any positive integer $k$ the map $F_*^{n \cdot 2^k}$ has a heteroclinic orbit $\cO_k$ between the periodic points  $\Lambda_*^{k}(p_0)$ and $\Lambda_*^{k-1}(p_0)$;
\end{itemize}
\end{theorem}

\begin{theorem}\label{MTHM2}
$G_* \equiv F_* \circ F_* \circ F_*$ admits a hyperbolic set $\cC_{G_*} \subset \cD_3$, where $\cD_3$ denotes the domain of $G_*$;
$$G_* \arrowvert_{\cC_G} \, { \approx \atop  \mbox{{\small \it  homeo}}   }  \, \sigma_2 \arrowvert_{\{0,1\}^{\fZ}},$$ 
whose Hausdorff dimension satisfies:
$$0.7673 \ge  {\rm dim}_H(\cC_{G_*}) \ge \varepsilon,$$
where $\varepsilon \approx 0.00013 \, e^{-7499}$ is strictly positive.  
\end{theorem}

\begin{figure}[h]
\begin{center}
\includegraphics[width=0.3\textwidth]{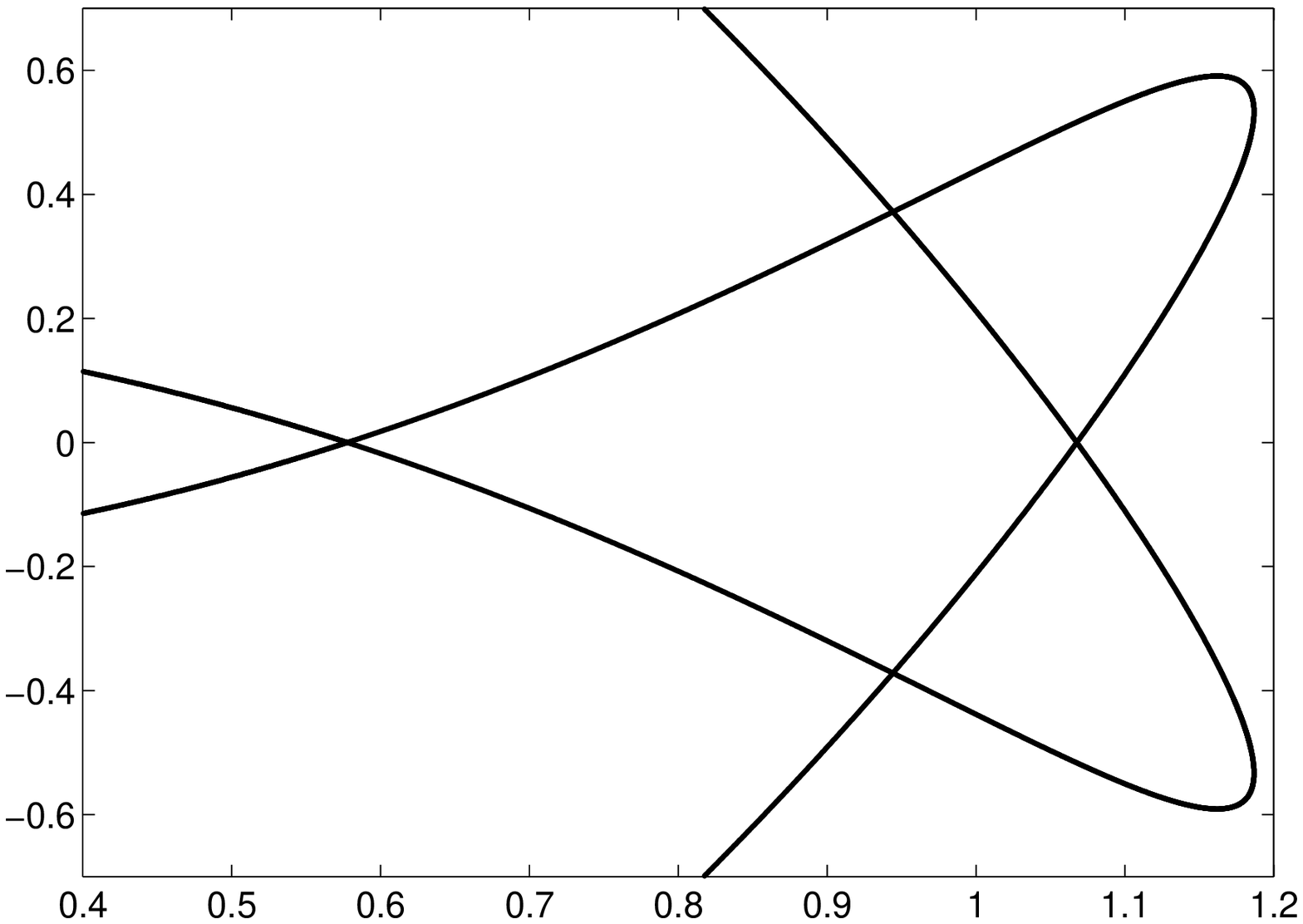}
\includegraphics[width=0.3\textwidth]{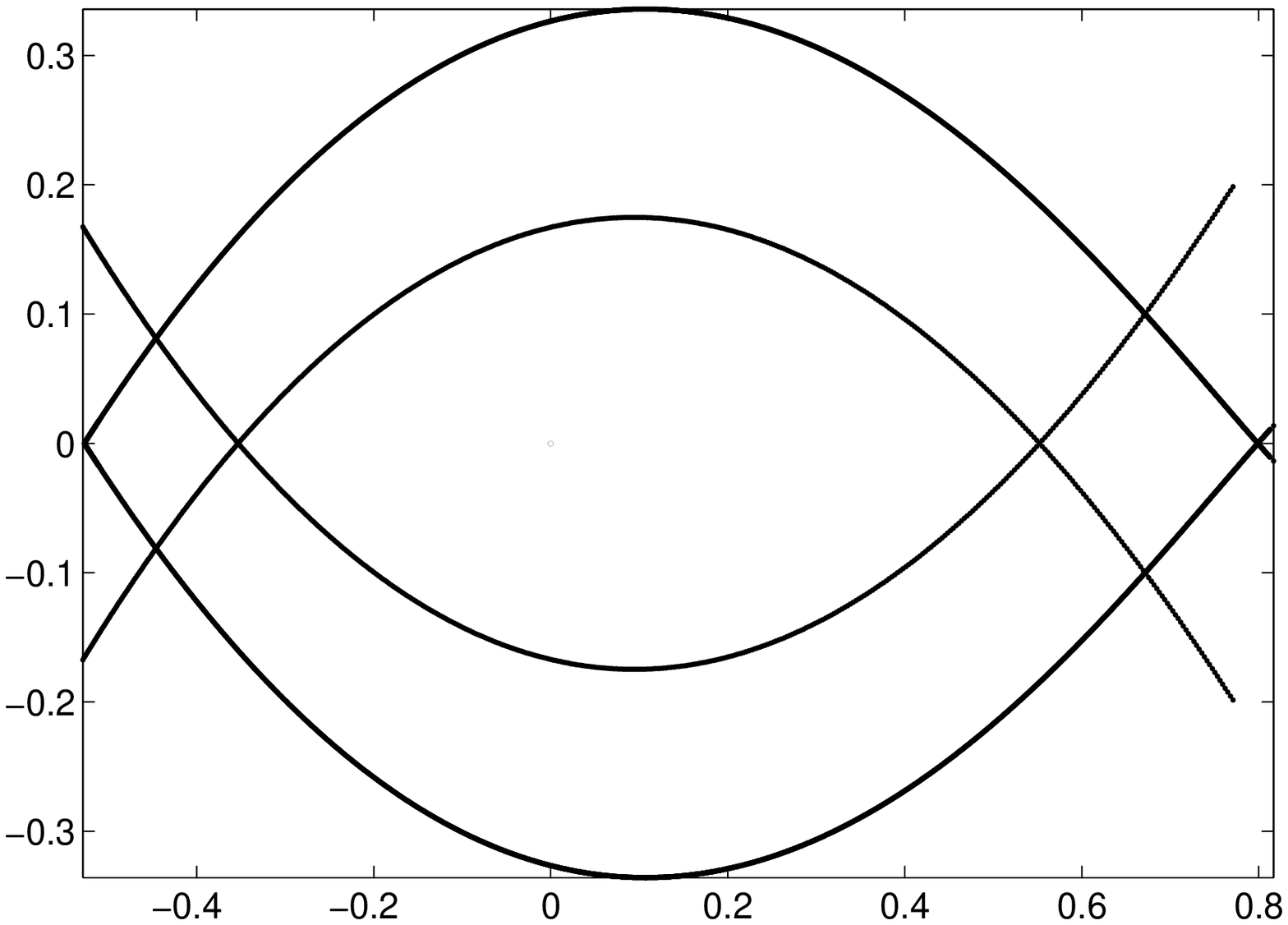} 
\caption{(a) The homoclinic tangle at the fixed point $p_0$. (b) The two components horseshoe of $G_*$.}\label{GHorsePic}\label{fHomoTraj}
\end{center}
\end{figure}

The transversal homoclinic intersection described in Theorem \ref{MTHM1} is illustrated in Figure \ref{fHomoTraj}(a), and the two component horseshoe described in Theorem \ref{MTHM2} is illustrated in Figure \ref{GHorsePic}(b). To summarize, in \cite{GJ1} we prove that the domain, $\cD$, of $F_*$ contains: a hyperbolic fixed point located at approximately $(0.577619,0)$, whose stable and unstable manifolds intersect transversally at a homoclinic point located at approximately $(1.067707,0)$; a two component horseshoe, one component contains the fixed point and the other component contains a period three point located at approximately $(-0.527155,0)$; a sequence of period $2^n$ points located at approximately $\lambda_*^n\times(0.577619,0)$, with transversal intersections of the stable and unstable manifolds of the points with period $2^{n-1}$ and $2^{n}$.

\section{Statement of results}
\begin{theorem}\label{MTHM}
The map $F_*$ has no elliptic islands of period less than $20$ in its domain $\tilde \cD$.
\end{theorem}

Most likely, the above result extends to all infinitely renormalizable maps, i.e., maps on the renormalization stable manifold, in a neighborhood of $F_*$. This would follow if one could prove that such maps are conjugated. Note, however, that the conjugacies between infinitely renormalizable maps constructed in \cite{GJ2} are only defined on the ``stable'' set, which is the closure of the orbit of $(0,0)$. It is not clear how, and if, these conjugacies can be extended to the entire invariant set. 

From the renormalization equation it follows immediately that if $x$ is a period $2n$ point, then $\Lambda^{-1} x$ is a period $n$ point. Hence Theorem \ref{MTHM} implies the following corollary.

\begin{corollary}\label{Cor1}
For $n=1,2,3,4$, the analytic continuation of the map $F_*$ to the domain $\cD^n$ has no elliptic islands of period less than $\lfloor 20/2^n \rfloor$ in $\tilde\cD^n$.
\end{corollary}

Formally, the opposite of Corollary \ref{Cor1} is also true: there are no elliptic islands of period $2^kn$, $k\in \mathbb{Z}_+$, $n\leq 20$ in $\Lambda^k \tilde\cD$. The invariant subset of $\Lambda^k \tilde \cD$ for any $k>0$, however, is empty, as can be proved using the methods from Section \ref{SM}. 

\begin{figure}[hp]
\begin{center}
\includegraphics[width=\textwidth]{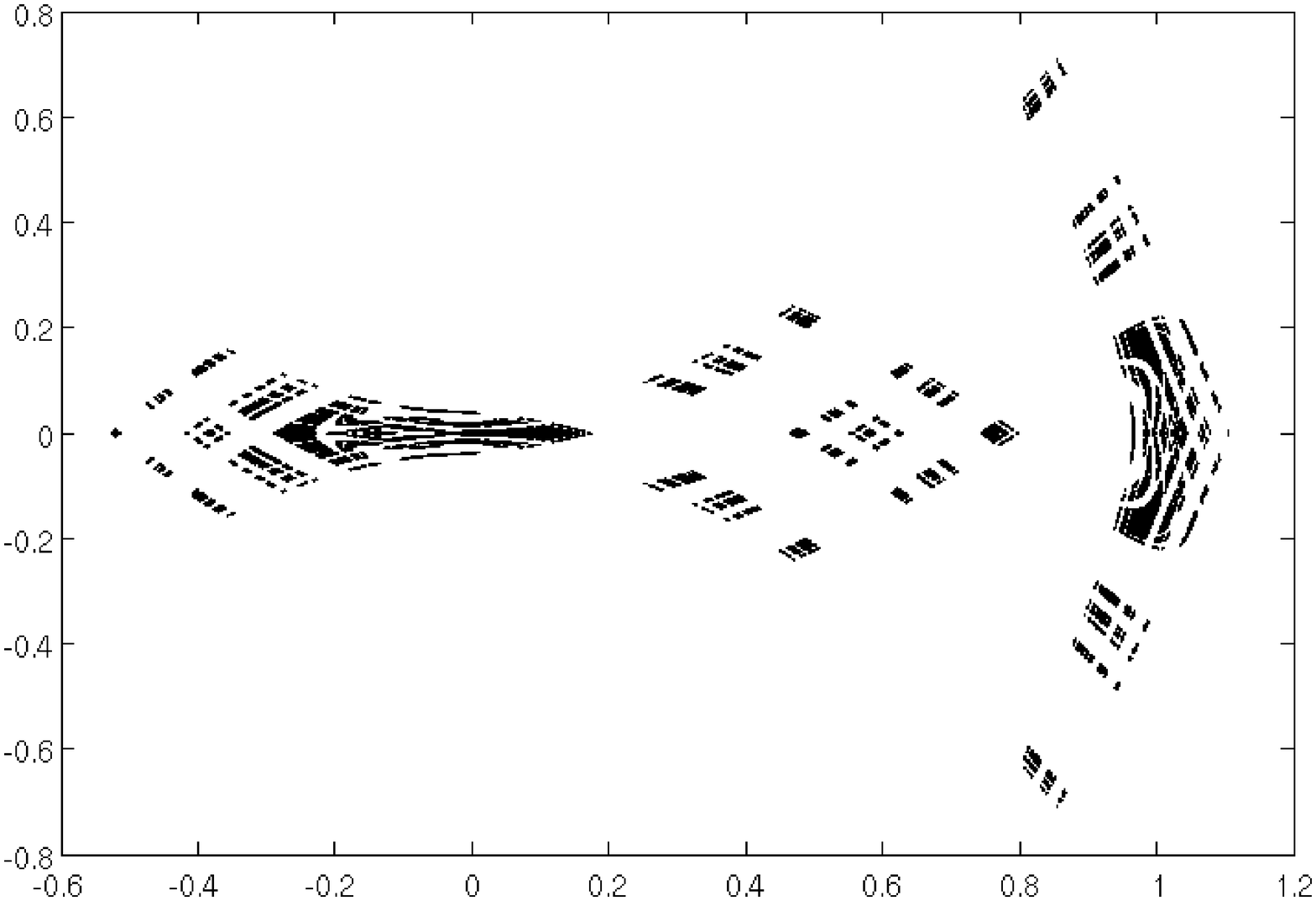}
\includegraphics[width=\textwidth]{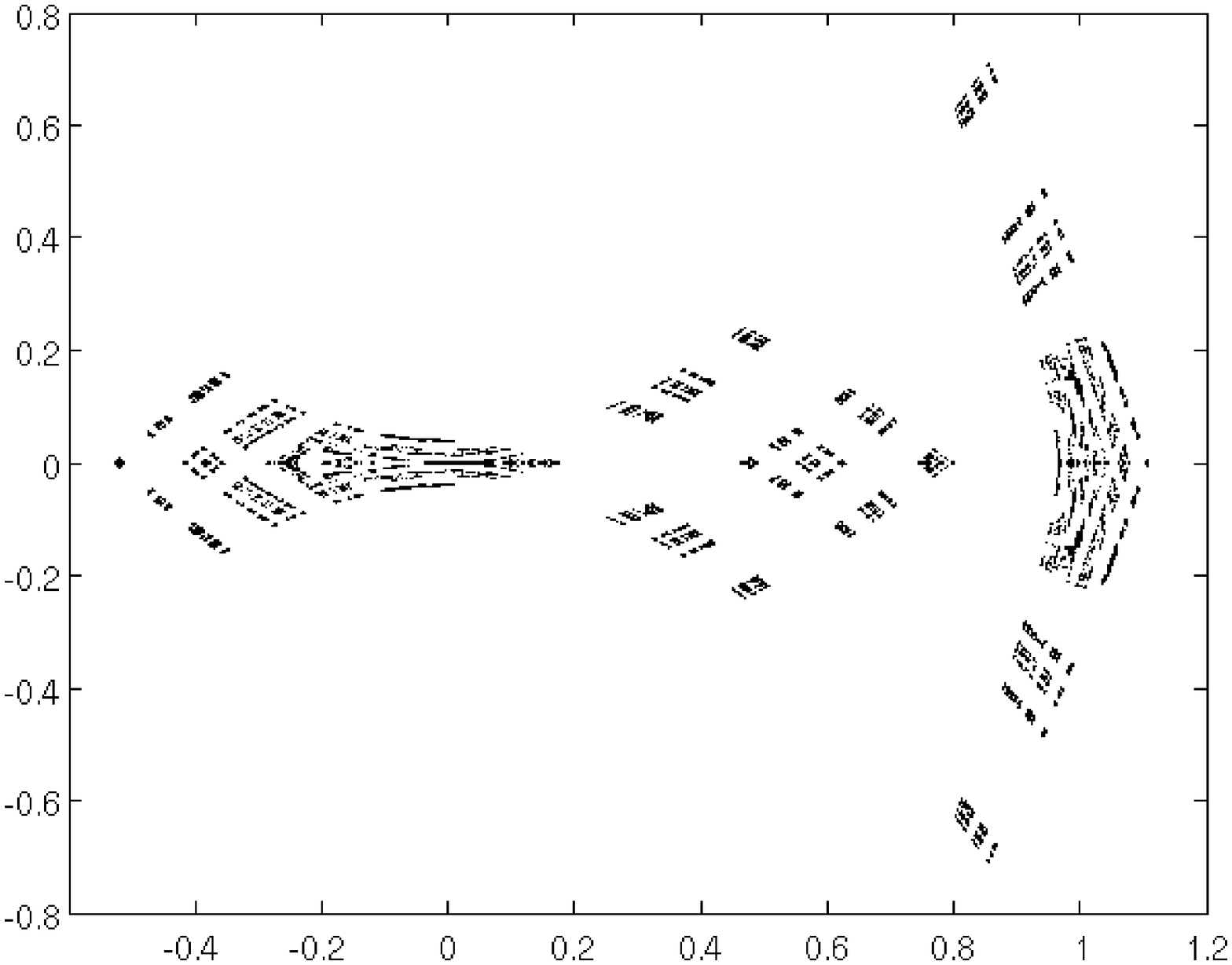} 
\caption{(a) The invariant subset of $\cD$, and (b) the periodic orbits up to order $19$. }\label{fInvDom}
\end{center}
\end{figure}

The invariant set, illustrated in Figure \ref{fInvDom}(a), consists of three regions, from left to right: the left component, i.e., the one containing the period three point, of the two component horseshoe; the sequence of rescaled periodic points with period $2^n$, and the associated multi component horseshoe induced by the transversal intersections of their stable and unstable manifolds; the topological horseshoe associated with the transversal homoclinic intersection. All of these regions are of a hyperbolic character, and contain, according to Theorem \ref{MTHM}, no elliptic islands of low period. In addition, the invariant set also contains the ``stable'' set from \cite{GJ2}.

To give further evidence that there are no elliptic islands, we have computed upper bounds on the measure of the invariant subset of $\cD$. As is seen from Figure \ref{fInvDom}(a), most of the remaining measure of the cover is around the sequence of rescaled periodic points, and their images. We note that the image of the point $(0,0)$ is close to $(1,0)$.

\begin{theorem}\label{TM}
The Lebesgue measure of the invariant subset of $\tilde \cD$ is less than $1.4594\%$ of the measure of $\tilde \cD$.
\end{theorem}

We will describe in Section \ref{SCompRes} that as the measure of the invariant set is estimated with higher resolution the estimate scales uniformly with the discretization size, this indicates that the measure should go to zero with the discretization size. It also shows that, up to the discretization size there are no elliptic islands, since once the discretization size is smaller than an elliptic island, the estimate of the measure would not decrease at the same rate. We therefore make the following conjecture, a similar conclusion for the renormalization limit of quadratic area-preserving maps was made in \cite{M82,M93}. Note that this statement is saying that the set of escaping orbits has full measure in the domain.

\begin{conjecture}\label{CM}
The Lebesgue measure of the invariant subset of $\tilde\cD$ is $0$.
\end{conjecture}

Based on Theorems \ref{MTHM} and \ref{TM}, and the structure of the invariant set, compared to the structure of the periodic orbits, see Figure \ref{fInvDom} (a) and (b), we make the following conjecture.  

\begin{conjecture}\label{MC}
The analytic continuation of the map $F_*$ to the domain $\cD^\infty$ has no elliptic islands in $\tilde\cD^\infty$.
\end{conjecture}

\section{Method}\label{SM}
In this section we describe the algorithms, and their mathematical justifications, that are used to prove Theorems \ref{MTHM} and \ref{TM}. Our proof uses interval arithmetic, see e.g. \cite{Ne90}, to rigorously compute enclosures of the action of the map $F_*$. Based on these enclosures we can construct a combinatorial representation of the map using graph techniques, i.e., periodic points of the map correspond to cycles in the graph of the combinatorial function. This approach is standard for dissipative systems and used in e.g. \cite{G01,G02a,G02b,TW09}. Since we study a conservative system there are additional difficulties, primarily the fact that once a low order cycle is detected, it is impossible to find a trapping region for it. Such information is used in \cite{G01,G02a,G02b,TW09} to adaptively reduce the number of nodes in the graph. In the conservative case we have to keep all nodes in the graph for all lengths of the cycles, there is nothing that prevents the map from having a higher order periodic orbit in a previously computed lower order cycle. This significantly increases the complexity of the algorithm. Furthermore, the essentially one dimensional discretization method that we use, does not lend itself to a nonuniform adaptive partition of the domain, as suggested in e.g. \cite{G02b}, instead we use a uniform adaptive discretization of the domain. A computational novelty of our paper is the use of generating functions in this setting.

Our algorithm consists of three steps, which are described in the subsections below. First we describe how the reduction to a combinatorial function is performed, second we describe how the candidate cycles are contracted, and in the final subsection we describe how we check for hyperbolicity. The algorithms described in this section are general, and can be used to prove the nonexistence of elliptic islands, up to some period, for any area-preserving map. The results of the computation for the universal area-preserving map associated with period doubling are given in the next section.

An implementation \cite{JP} of the algorithm described below has been made using the CAPD library \cite{CAPD} for interval arithmetic operations.

\subsection{Combinatorial representation}
The idea behind the combinatorial representation is simple, partition the domain, compute the range of the map restricted to each set in the partition, and check which sets in the partition that the image intersects. Once all these ranges are computed, a graph is constructed in which each set in the partition is represented by a node, and each intersection is represented by an edge. 

Recall our definition of an area-preserving map, it is an exact symplectic reversible diffeomorphism of a subset of ${\fR}^2$ onto its image. Thus, the action of $F_*$ is described by a generating function. Therefore, it is enough to discretize one of the two dimensions. With a physical interpretation of the generating function, we only discretize the position variable, not the momentum.  

Consider three consecutive points, $(x,p_x), (y,p_y), (z, p_z)$, on an orbit of $F_*$, with position values $x,y,$ and $z$. The following holds, 
\begin{equation}\label{gen_func_per}
\left( x \atop -s(y,x) \right) 
{{ \mbox{{\small \it  F}} \atop \mapsto} \atop \phantom{\mbox{\tiny .}}}  
\left( y \atop s(x,y) \right) = \left( y \atop -s(z,y)\right) 
{{ \mbox{{\small \it  F}} \atop \mapsto} \atop \phantom{\mbox{\tiny .}}}  
\left( z \atop s(y,z) \right).
\end{equation}
Thus, for any three consecutive points on the trajectory the following equation is satisfied:
\begin{equation}\label{three_points}
 s(x,y)+s(z,y)=0.
\end{equation}
Given a pair of points $(x,y)$; (\ref{three_points}) is a nonlinear equation in $z$ that can be solved by e.g. the interval Newton method.

To build the graph, we first uniformly discretize the projection of the domain on the first variable. This is an interval, which we denote by $I$, the discretization is denoted by $U=\{I_i\}_{i=1}^N$. The combinatorial representation is a $N\times N$-graph, where an edge between $(I_i,I_j)$ and $(I_j,I_k)$ represents that there exists a trajectory whose projection onto the first coordinate visits $I_i$, $I_j$, and $I_k$, in that order. For each pair $(I_i,I_j)$, (\ref{three_points}) is solved using the interval Newton method. For each set $I_k$ with a non-empty intersection with the result, an edge is added to the graph. 

From the graph, all cycles up to some given length are computed by an exhaustive depth first search algorithm. All periodic points of the map are contained in a cycle whose length divides the period. Due to overestimation and discretization errors, however, most cycles are spurious and do not actually contain periodic orbits. All cycles detected in the graph are investigated further, as described in the next subsection.

To decrease the computing time, the search for cycles in the graph is preceded by a contraction step, similar to the one in \cite{G02a,G02b}, where all nodes that do not have both inward and outward pointing edges are removed. This procedure is performed since an edge that goes to a node which leaves the domain cannot be part of an orbit between points in the invariant set. Similarly, a node that is unreachable is not in the invariant set. Each time a node is removed, all edges pointing to or from it are also removed. This procedure is iterated until no more nodes are removed. The remaining nodes represent the invariant set of the domain. 
In our case, this reduction step removes about $90$\% of the nodes and edges from the graph.

\subsection{Contracting the candidate cycles}
Given a period $n$, the first step of the algorithm produces a set $C_n$ of $N_n$ sequences of $n$ intervals from the partition 
\begin{equation*}
C_n := \{I^i\}_{i=1}^{N_n} = \{(I^i_1,...,I^i_n) \}_{i=1}^{N_n}. 
\end{equation*}
Each interval $I^i_j\in U$, and each sequence $(I^i_1,...,I^i_n)$, corresponds to a cycle $\left((I^i_1,I^i_2),(I^i_2,I^i_3),...,(I^i_{n-1},I^i_n),(I^i_n,I^i_1)\right)$ in the graph constructed in the previous section. That is, the set $C_n$ is the set of all sequences, $I^i$, of intervals, such that there exists a cycle in the graph, whose nodes correspond to neighboring pairs of intervals in $I^i$, in the same order.

To quickly remove most of the spurious cycles we contract each sequence by solving $n$ one dimensional problems, constructed by solving for $y$ in (\ref{three_points}), using the interval Newton method. Below $I=I^{i}$, for some $i$. I.e., we iterate: 
$$
I_i=I_i\cap\left(\textrm{mid}(I_i)-(s_2(I_{i-1},I_i)+s_2(I_{i+1},I_i))^{-1}(s(I_{i-1},\textrm{mid}(I_i))+s(I_{i+1},\textrm{mid}(I_i)))\right),
$$
where $\textrm{mid}(I_i)$ denotes the midpoint of the interval $I_i$. If any intersection in a cycle is empty, the entire corresponding cycle is removed. This step removes the vast majority of candidate cycles. The remaining cycles are contracted more tightly on the candidate periodic orbit using the interval Krawczyk method. We define:
\begin{equation}\label{fZero}
Z(I_1,\dots,I_n)=\left[\begin{array}{c}
s(I_n,I_1)+s(I_2,I_1) \\
\cdots \\
s(I_{n-1},I_n)+s(I_1,I_n)
\end{array}\right],
\end{equation}
its derivative is given by:
\begin{tiny}
\begin{eqnarray}\label{dfZero}
\hspace{-2cm}
DZ(I_1,\dots,I_n) = \\ \nonumber \hspace{-1cm} {\left[\begin{array}{cccccc}
s_2(I_n,I_1)+s_2(I_2,I_1) & s_1(I_2,I_1) & 0 & \cdots & 0 & s_1(I_n,I_1) \\
s_1(I_1,I_2) & s_2(I_1,I_2)+s_2(I_3,I_2) & s_1(I_3,I_2) & \cdots & 0 & 0\\
& & \cdots & & & \\
s_1(I_1,I_n) & 0 & 0 & \cdots & s_1(I_{n-1},I_n) & s_2(I_{n-1},I_n)+s_2(I_1,I_n)
\end{array}\right]}
\end{eqnarray}
\end{tiny} 
Zeros of (\ref{fZero}) correspond to periodic cycles of $F_*$. Note that (\ref{dfZero}) is a symmetric matrix, since $s_1$ is a symmetric function. The interval Krawczyk operator is defined as:
$$
K(I) := \textrm{mid}(I)-CZ(\textrm{mid}(I))+(Id-CDZ(I))(I-\textrm{mid}(I)),
$$ 
where $C$ is any nonsingular matrix; if possible we use $C=(DZ(\textrm{mid}(I)))^{-1}$. Similarly as with the Newton operator above, we contract the cycles by iterating the following intersection:
$$
I=I\cap K(I)
$$
If the intersection of a cycle is empty, it is removed. See \cite{Ne90} for a thorough description of the interval Newton and Krawczyk methods. 

Note that we do not check for existence of periodic orbits, even though the Krawczyk method can be used for this. The reason is that we are only interested in proving that there are no elliptic cycles, so if we additionally prove that some non-periodic length $n$ orbits are hyperbolic that is irrelevant. Furthermore, given the discretization size, a candidate cycle might contain several periodic orbits of the same length, and to decrease the complexity of the algorithm we want to avoid having to split this interval in several subintervals, which would be required if we wanted to prove existence of periodic orbits. 

\subsection{Proving hyperbolicity}
Given that the maps that we study are area-preserving, the only possible eigenvalues of the derivative at a periodic point are real, or complex and of unit modulus. We compute an enclosure of the derivative along each candidate cycle using equation (\ref{Fder}), and check that the eigenvalues are real and not equal. In the case when we are unable to prove whether a cycle is hyperbolic or elliptic, we split the widest enclosure, $\tilde I_i$, of a point in the orbit. Note that, $\tilde I_i\subset I_i$, for some $I_i$ in the cycle, since the original partition is contracted.

\section{Computational results}\label{SCompRes}
In this section we describe the results of the computations using the algorithm presented in the previous section on the universal area-preserving map associated with period doubling, $F_*$. These computations constitute the proofs of Theorems \ref{MTHM} and \ref{TM}. Computations were made on a 32 CPU computer (8 2.2 GHz Athlon Quad Processors) with 32Gb of RAM, and on a 3.2 GHz Dual-Core AMD Opteron Processor with 66Gb of RAM. The C++ code of the implementation is available for download at \cite{JP}. 

\subsection{The invariant set}
The projection of the domain, $\tilde \cD$, on the first coordinate is $I=[-1.1,2.1]$. We start with a coarse discretization, using 100 subintervals, in order to find out whether the invariant set is all of $I$ or only a subset. We repeat this procedure a few times until the results stabilize. This leaves us with the invariant interval $[-0.55,1.13]$, which is the interval that we will study for the remainder of this paper. We set 
$$
I=[-0.55,1.13].
$$

To estimate the measure of the invariant set, we compute successive covers of it, as described in the previous section. We start with the discretization size $400$. Each node, $(I_i,I_j)$, in the graph corresponds to a set, $(I_i,-s(I_j,I_i))$, in the cover of the invariant set. We estimate the area of each box, and sum over all nodes in the graph to get the upper estimate of the cover. Given a discretization size, we split each node in two, regenerate the edges, and compute the new enclosure of the invariant set. The measure of the enclosures at different (nominal) discretization sizes are given in Table \ref{tM}. Each time the size of the sides of the covering boxes is decreased by a factor of two, the remaining area is decreased by a factor $1.5$. Such a uniform scaling is consistent with a cover of a Cantor set, giving strong evidence that there are no elliptic islands in the domain. These computations also shows that there are no elliptic islands in the domain that contain a ball with diameter larger than $0.00013125$, the final discretization size. Since, if such islands existed, the improvement would decrease, once this size was passed. Together with the bound $3.312$ on the measure of $\tilde\cD$ from Section 2, these computations finish the proof of Theorem \ref{TM}.

\begin{table}[h]
\begin{center}
\begin{tabular}{c|rrc}
Discretization (nominal) & Nodes & Enclosed Area & Improvement  \\ \hline
$400$ & $9366$ & $0.3742$ & $-$ \\
$800$ & $25051$ & $0.2492$ & $1.50$ \\
$1600$ & $67277$ & $0.1671$ & $1.49$ \\
$3200$ & $176304$ & $0.1097$ & $1.52$ \\
$6400$ & $459828$ & $0.07184$ & $1.53$ \\
$12800$ & $1227236$ & $0.04834$ & $1.49$
\end{tabular}
\caption{The upper estimates of the measure of the invariant subset of $\tilde \cD$. The nominal discretization size is the number of intervals with uniform size that the interval $[-0.55,1.13]$ would be split into if the discretization was non-adaptive, i.e., the first coordinate of each box in the cover has length $1.68/\tt{discretization}$. The last column shows the improvement of the estimate compared to the previous level.}
\end{center}
\end{table}\label{tM}

\subsection{Periodic orbits}
Before we start to compute enclosures of the periodic orbits we have to choose the (nominal) discretization size of $I$ that will be used for the main computations. After some trial and error, it turns out that the most efficient is to choose a discretization which is as fine as possible, since this reduces the number of spurious cycles and the number of cycles that must be split during the contraction part of the algorithm. Partitioning $I$ into $N$ subintervals means that we will study a (nominally) $N^2$-graph, so there is a memory bound on how large we can take $N$. The largest (nominal) $N$ that we are able to handle conveniently is $N=3200$, which is what we use in our proof. I.e. we set 
$$
I_i=[-0.55+0.000525i,-0.55+0.000525(i+1)], \quad \textrm{ for }i=0,...,3199.
$$ 
This means that, apriori, our graph has $10240000$ nodes and up to $32768000000$ edges. Fortunately, our combinatorial function is much smaller. Running the reduction procedure, which reduces the graph to the invariant set, at the finest resolution, yields a graph with $117883$ nodes and $413556$ edges. The corresponding invariant set is illustrated in Figure \ref{fInvDom}(a). Ideally, since we study a map, as the discretization size goes to zero, the number of edges per node should go to one. With the discretization at hand, we get about $3.51$ edges per node.

For each period we locate all the corresponding cycles, contract them, and check for hyperbolicity as explained in the previous section. The memory bound part of the algorithm is a pre-calculation reduction step that locates which nodes that can be part of a cycle of a given length (period $19$ requires approximately $100$Gb). We only search for cycles in these reduced graphs. Each reduced graph is much smaller than the graph itself, the largest subgraph is the one for period $18$, which has $25631$ nodes and $69956$ edges. The number of nodes and edges of each subgraph, and the number of cycles, $b_n$, that we found in the graph are given in table \ref{tCycles}, together with the number of cycles, $h_n$, that remained after the contraction procedures. We only checked for hyperbolicity on the cycles that remained. The actual number of cycles that we had to compute was between $b_n$ and $b_n/n$, since a cycle of length $k$ appears $k$ times, but we only compute it once. This is done by keeping the nodes sorted, and traversing the graph so that no nodes with lower index than the starting index is visited. This also implies that we can start searching for cycles from all nodes in parallel, which significantly speeds up the algorithm.

The cycles up to order $19$ that we checked for hyperbolicity are illustrated in Figure \ref{fInvDom}(b). Note that we do not claim that all of these cycles contain hyperbolic periodic points; what we prove is that all periodic points of order less than $20$ are contained in these cycles, and that all periodic orbits up to order $19$ that are contained in those cycles are hyperbolic. This concludes the proof of Theorem \ref{MTHM}.

\begin{table}[h]
\begin{center}
\begin{tabular}{c|rrrrr}
$n$ & $\#$nodes & $\#$edges & $b_n$ & $\left\lfloor\frac{b_n}{n}\right\rfloor$ & $h_n$ \\ \hline
$1$ & $1$ & $1$ & $1$ & $1$ & $1$ \\ 
$2$ & $15$ & $27$ & $15$ & $7$ & $2$ \\
$3$ & $15$ & $28$ & $22$ & $7$ & $3$ \\
$4$ & $146$ & $425$ & $419$ & $104$ & $5$ \\
$5$ & $52$ & $92$ & $206$ & $41$ & $3$ \\
$6$ & $167$ & $385$ & $1938$ & $323$ & $8$ \\
$7$ & $117$ & $234$ & $1807$ & $258$ & $3$ \\
$8$ & $1663$ & $5336$ & $212547$ & $26568$ & $29$ \\
$9$ & $2135$ & $5811$ & $265342$ & $29482$ & $16$ \\
$10$ & $787$ & $1838$ & $188230$ & $18823$ & $14$\\
$11$ & $1060$ & $2270$ & $377323$ & $34302$ & $21$ \\
$12$ & $5072$ & $14604$ & $15861950$ & $1321829$ & $187$ \\
$13$ & $1691$ & $3425$ & $2123239$ & $163326$ & $33$ \\
$14$ & $4275$ & $9893$ & $30629187$ & $2187799$ & $104$ \\  
$15$ & $5106$ & $11360$ & $28118492$ & $1874566$ & $273$ \\ 
$16$ & $12823$ & $35343$ & $4930683251$ & $308167703$ & $80901$ \\ 
$17$ & $11324$ & $25850$ & $253950540$ & $14938267$ & $220$ \\
$18$ & $25631$ & $69956$ & $24952575180$ & $1386254176$ & $13159$ \\ 
$19$ & $24558$ & $62806$ & $3638886917$ & $191520364$ & $435$
\end{tabular}
\caption{For each period $n$, the number of nodes and edges refer to the corresponding subgraph, $b_n$ is the number of cycles, including multiplicity, the actual number of cycles that are computed is between $b_n/n$ and $b_n$. $h_n$ is the number of cycles that remained after the contraction and were checked for hyperbolicity.}
\end{center}
\end{table}\label{tCycles}

\begin{figure}[hp]
\begin{center}
\includegraphics[width=0.325\textwidth]{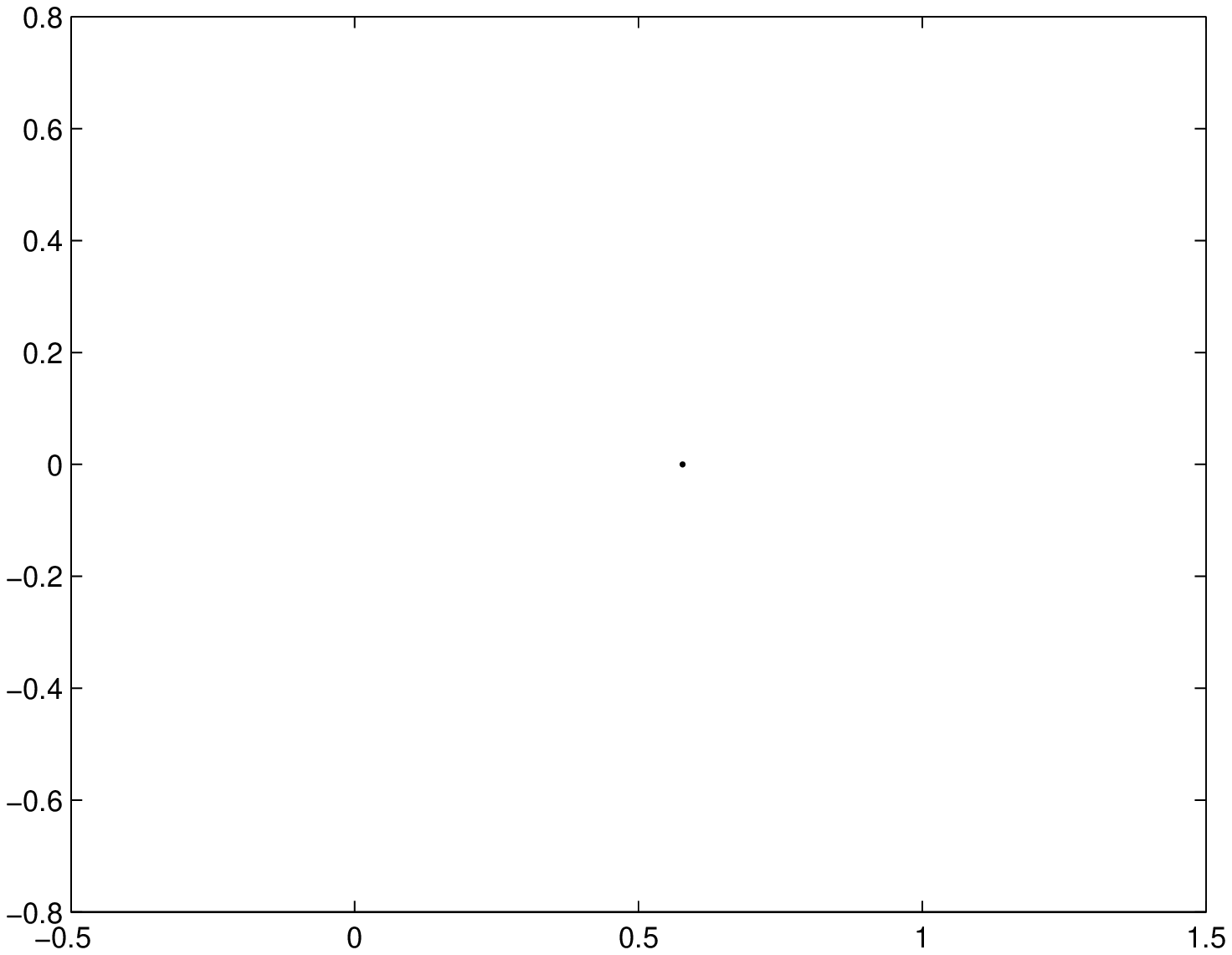}
\includegraphics[width=0.325\textwidth]{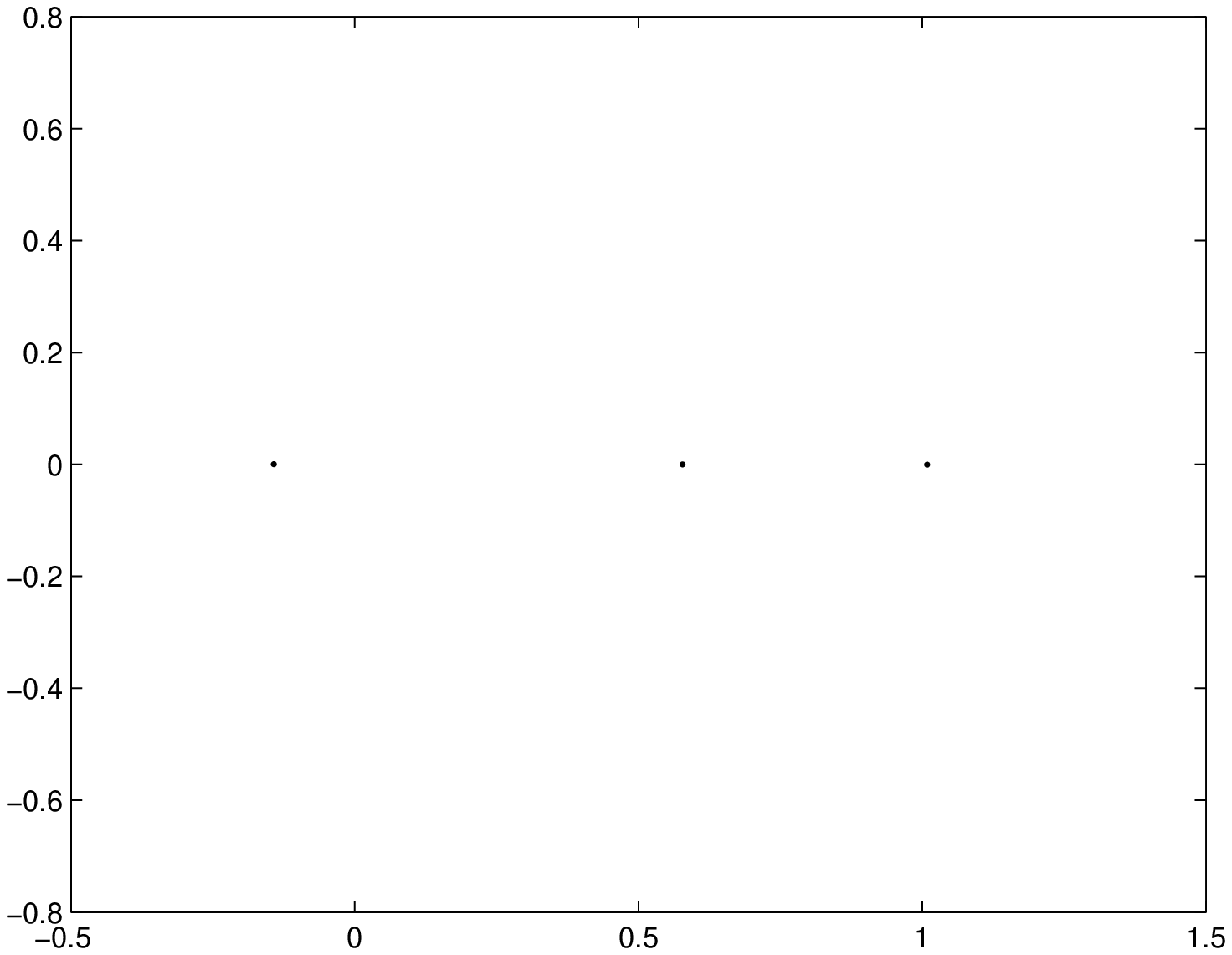}
\includegraphics[width=0.325\textwidth]{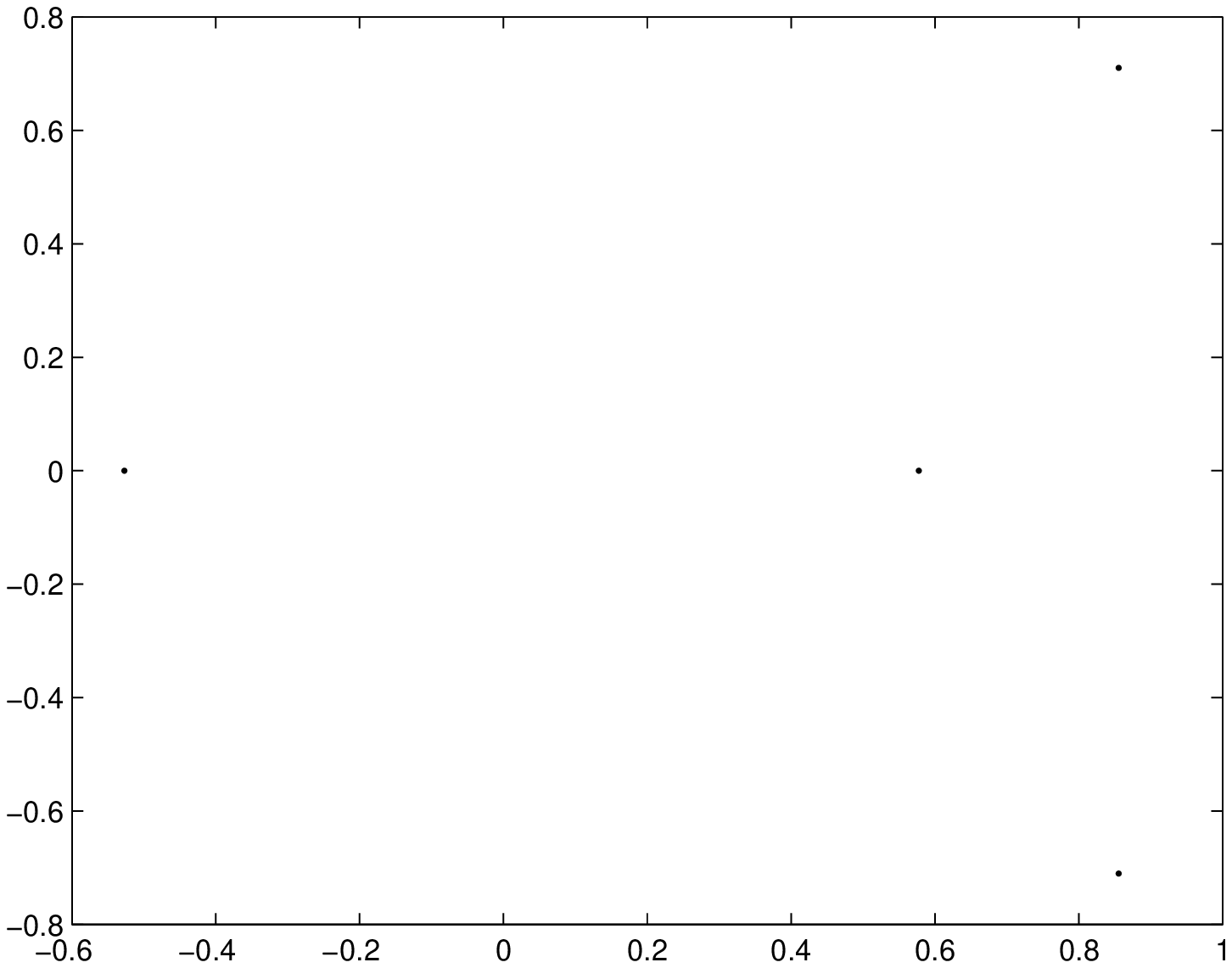}
\includegraphics[width=0.325\textwidth]{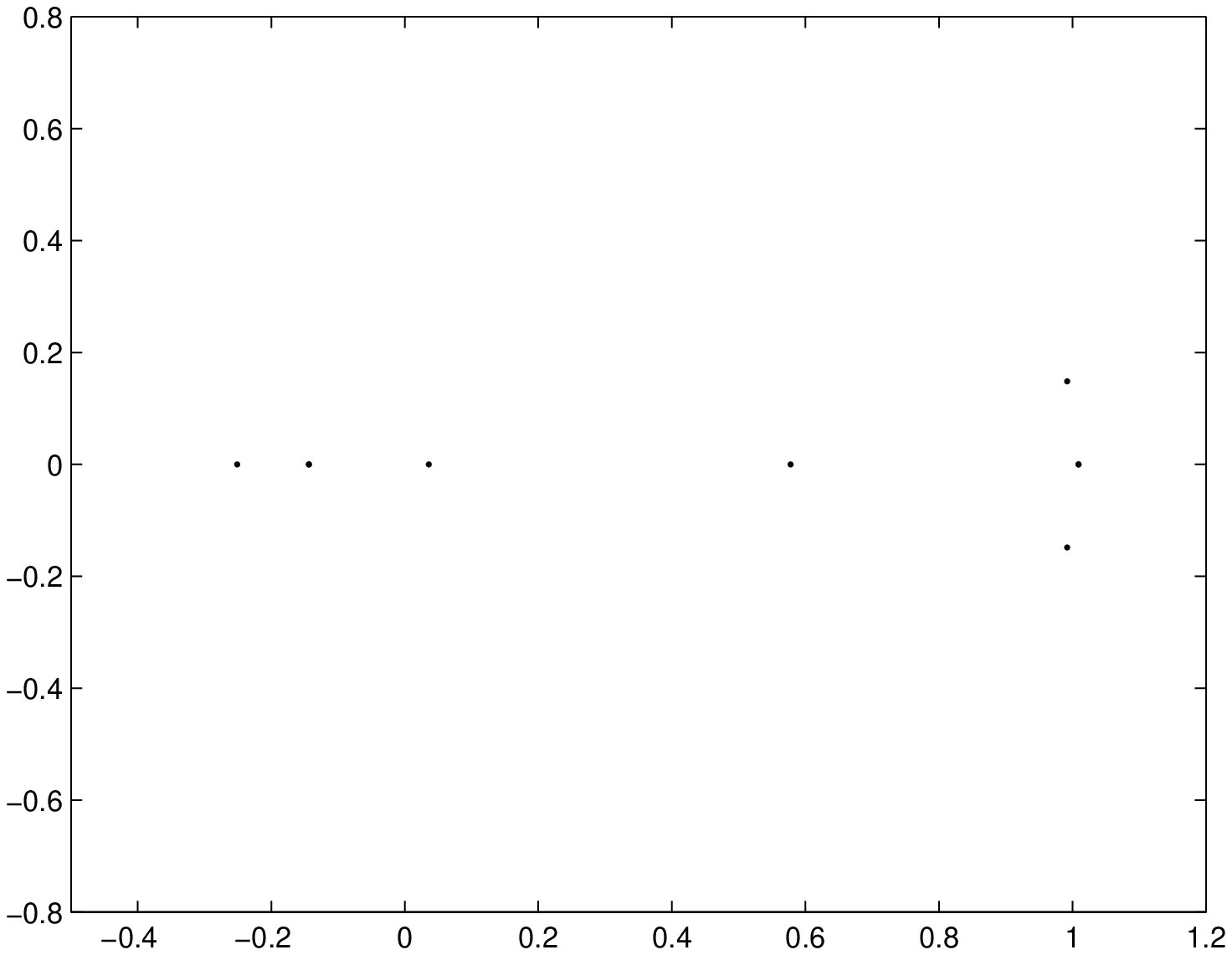}
\includegraphics[width=0.325\textwidth]{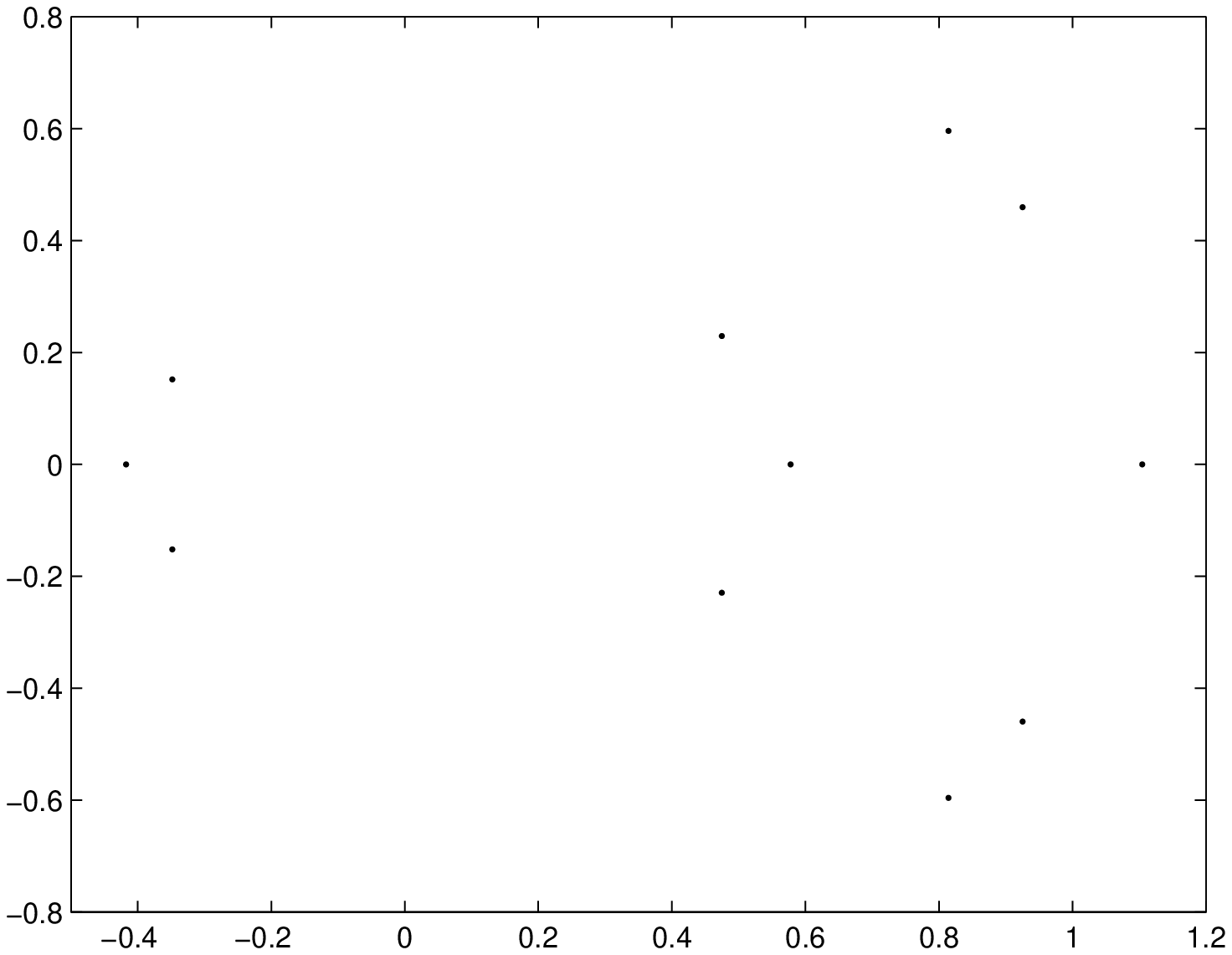}
\includegraphics[width=0.325\textwidth]{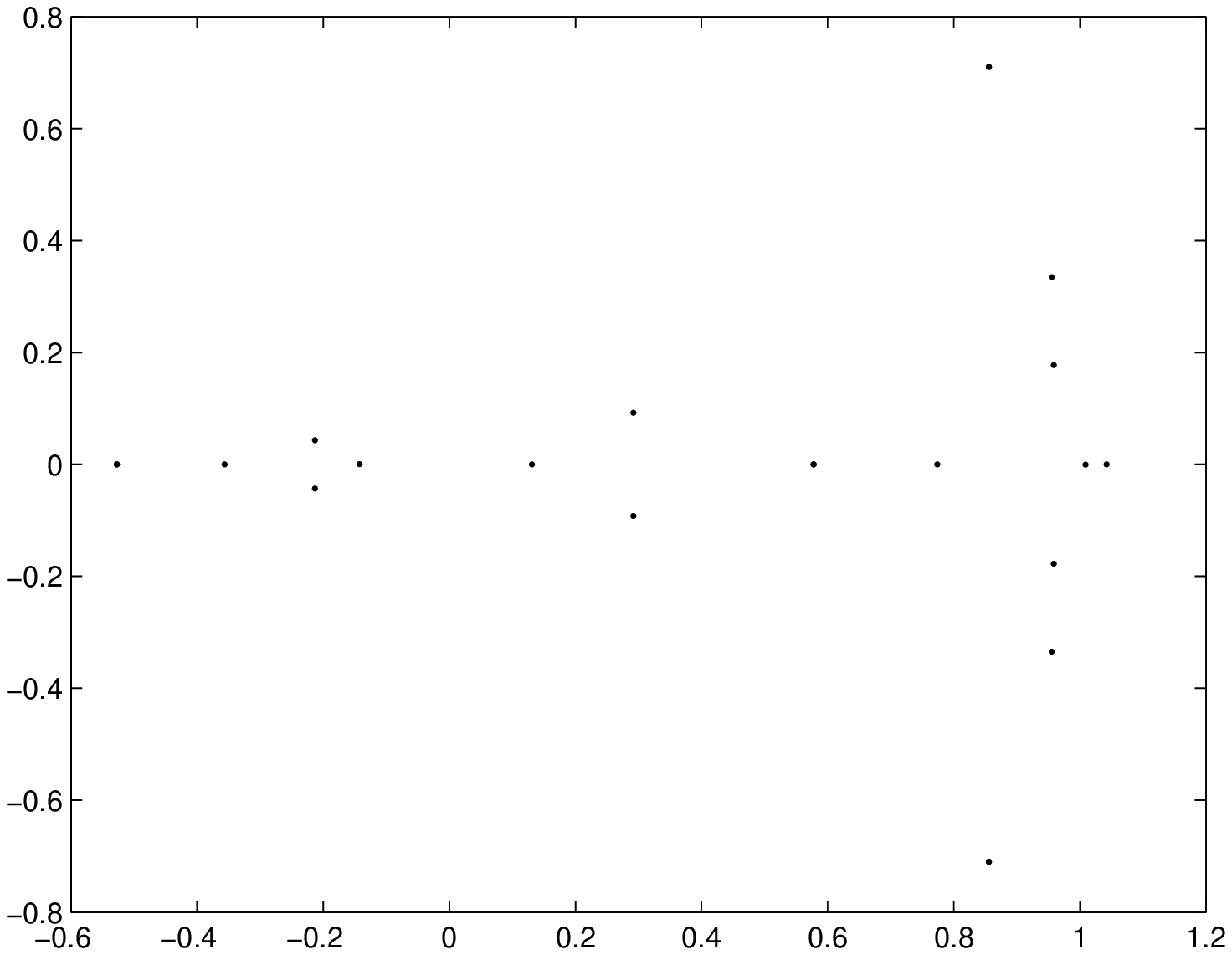}
\includegraphics[width=0.325\textwidth]{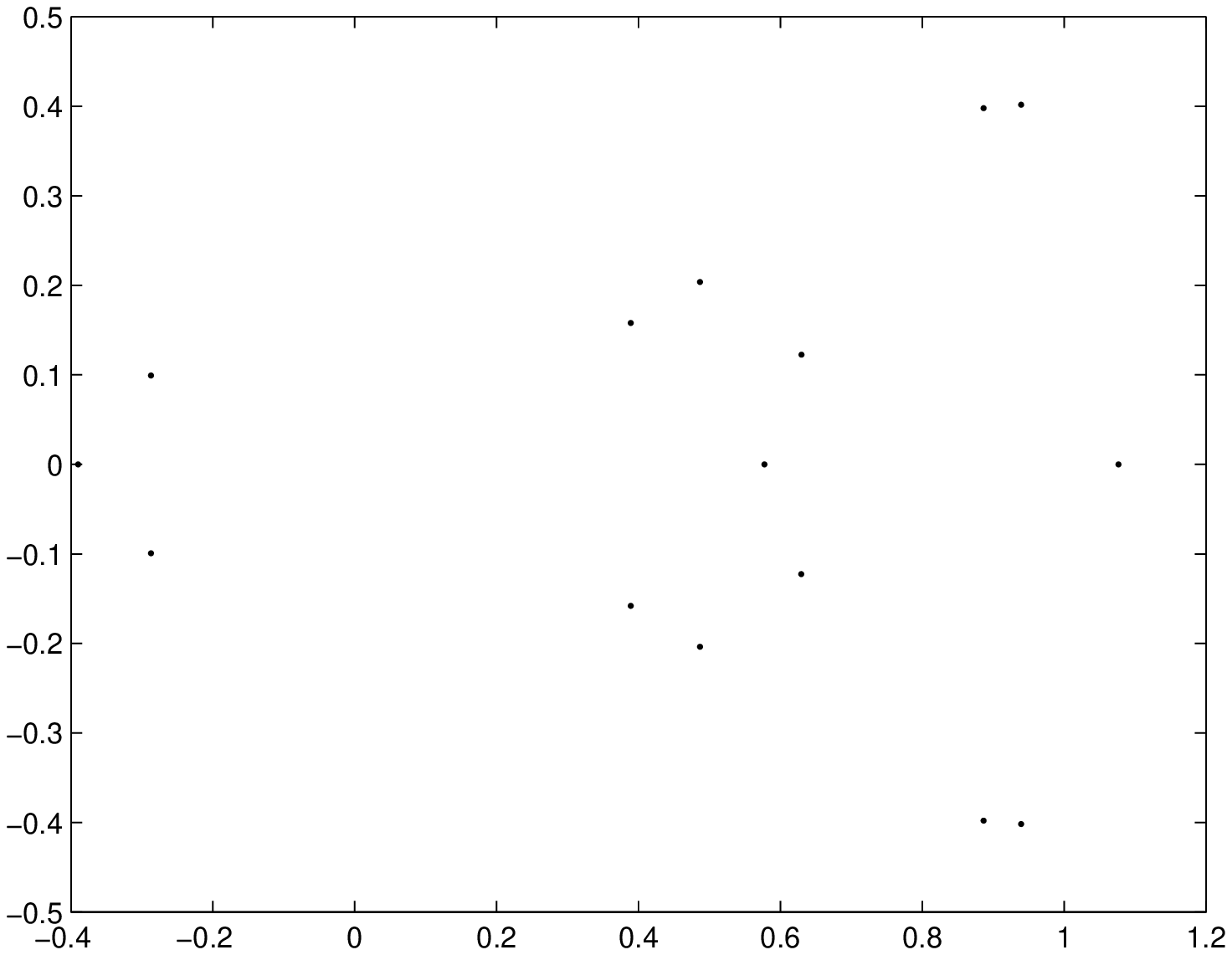}
\includegraphics[width=0.325\textwidth]{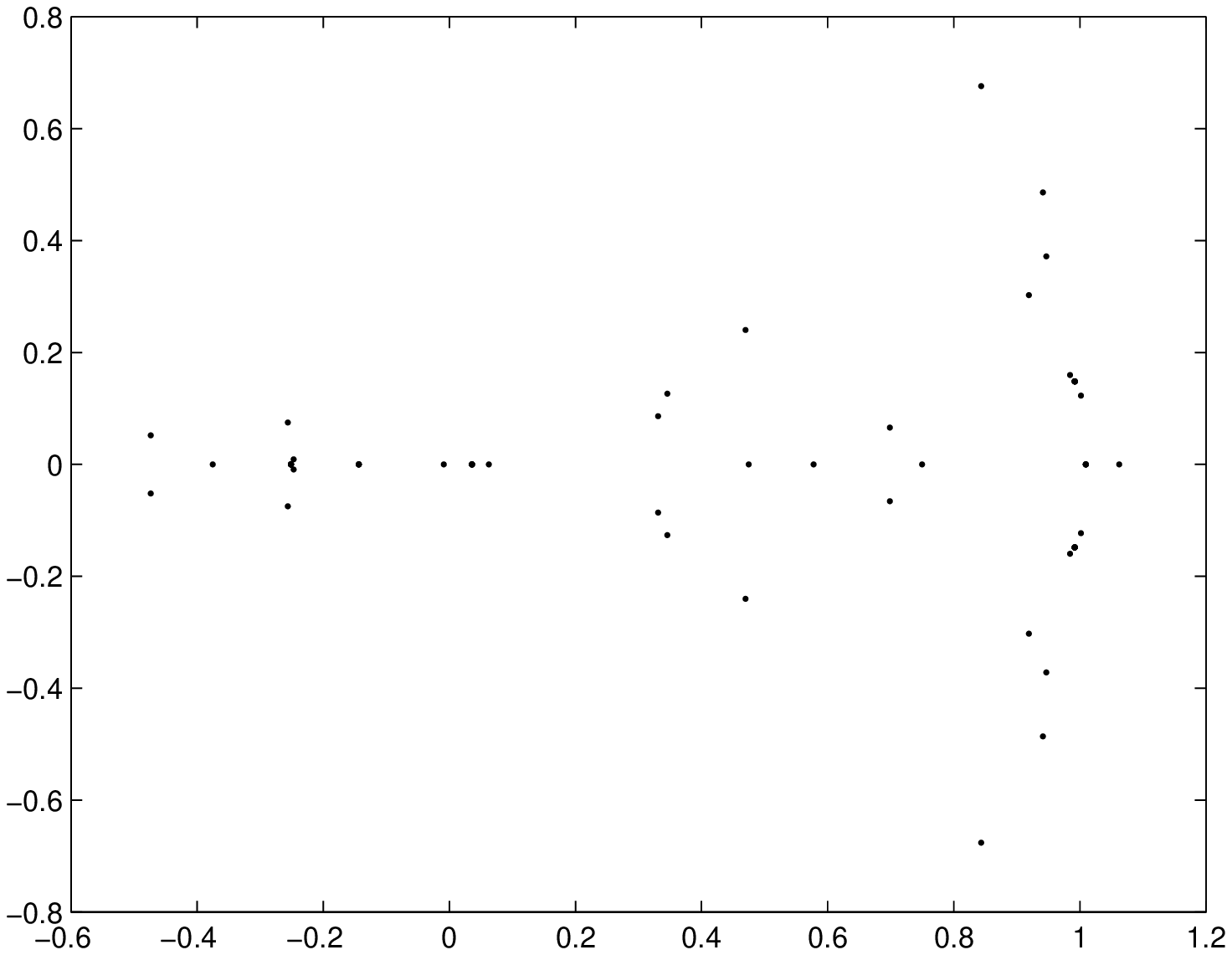}
\includegraphics[width=0.325\textwidth]{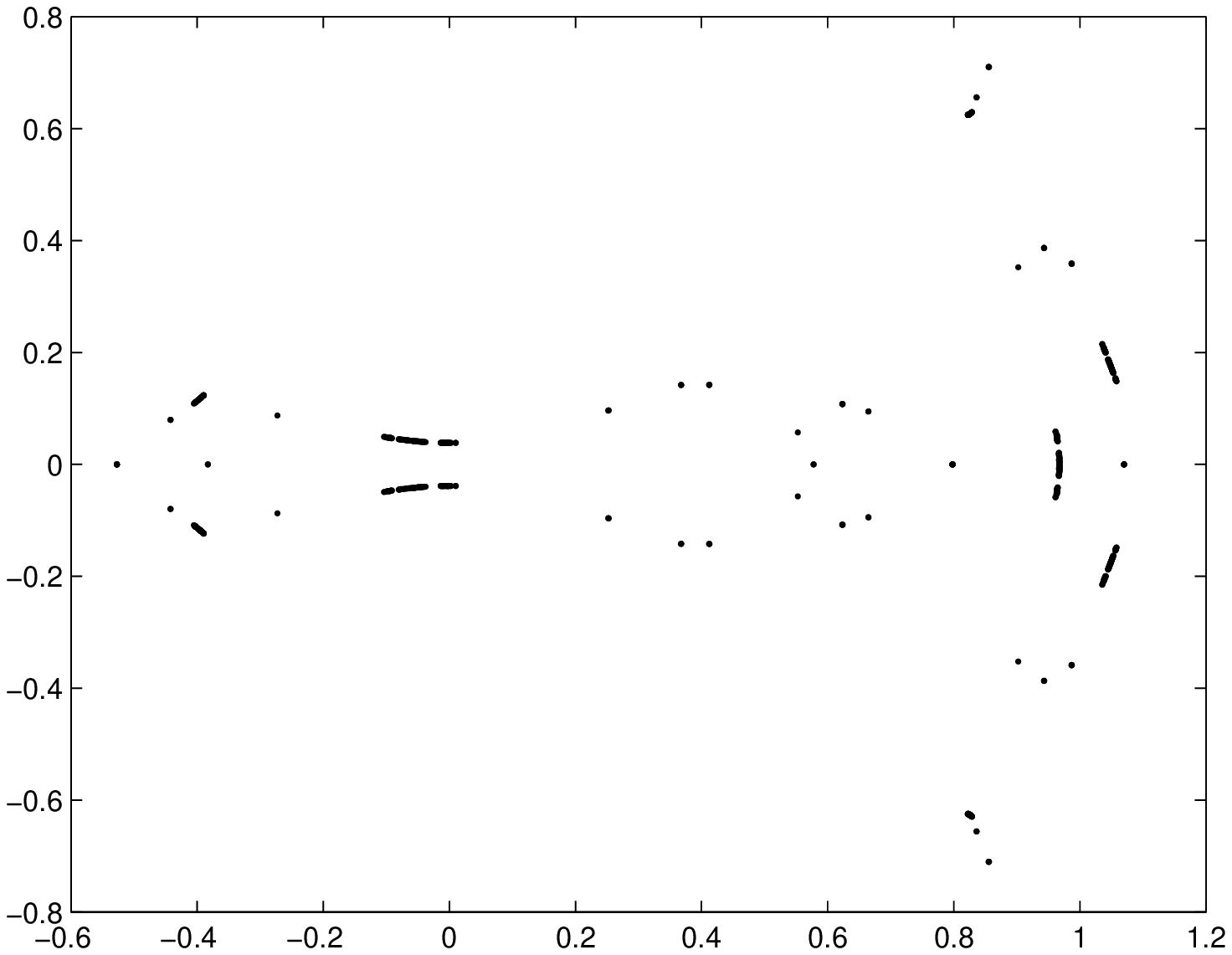}
\includegraphics[width=0.325\textwidth]{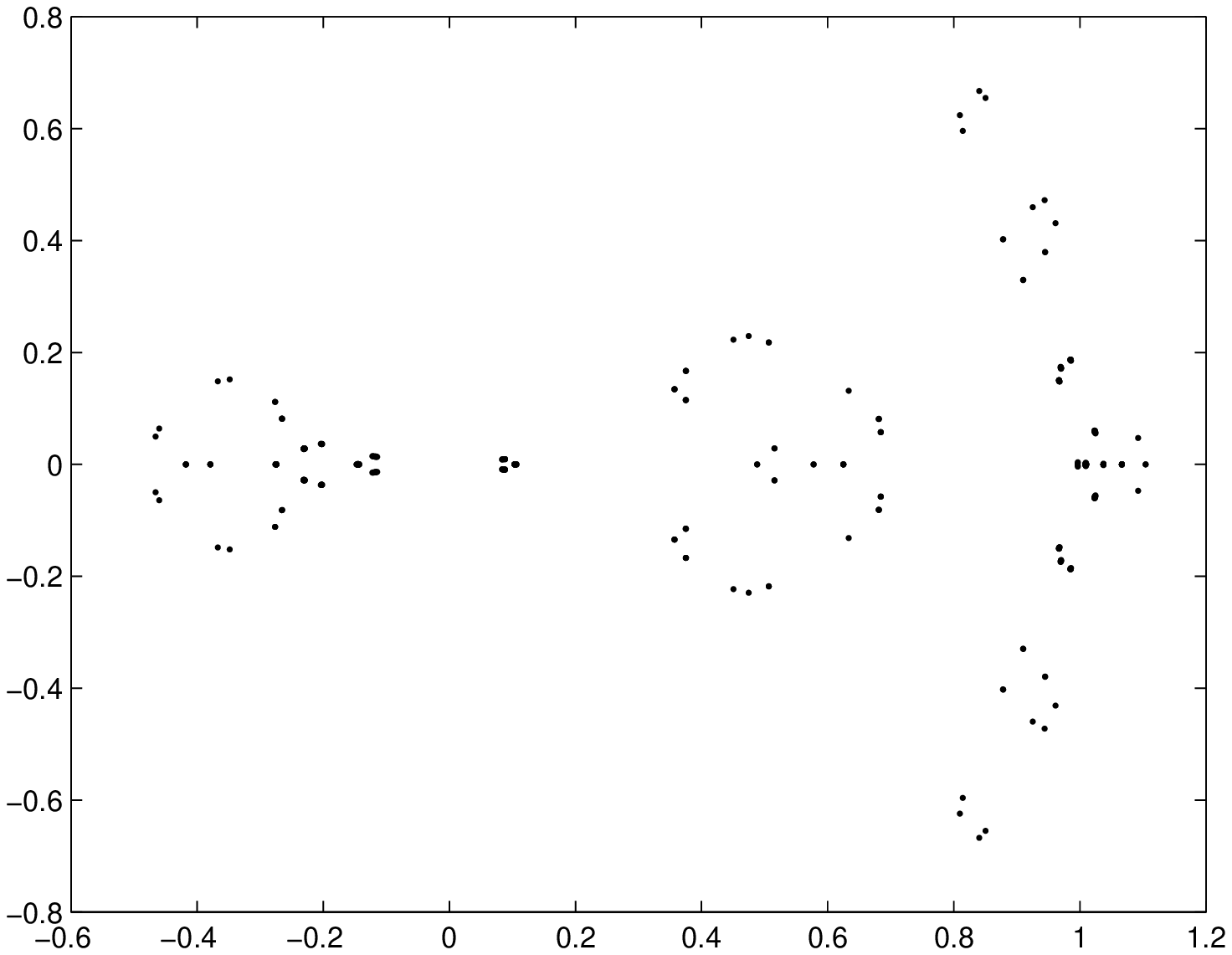}
\includegraphics[width=0.325\textwidth]{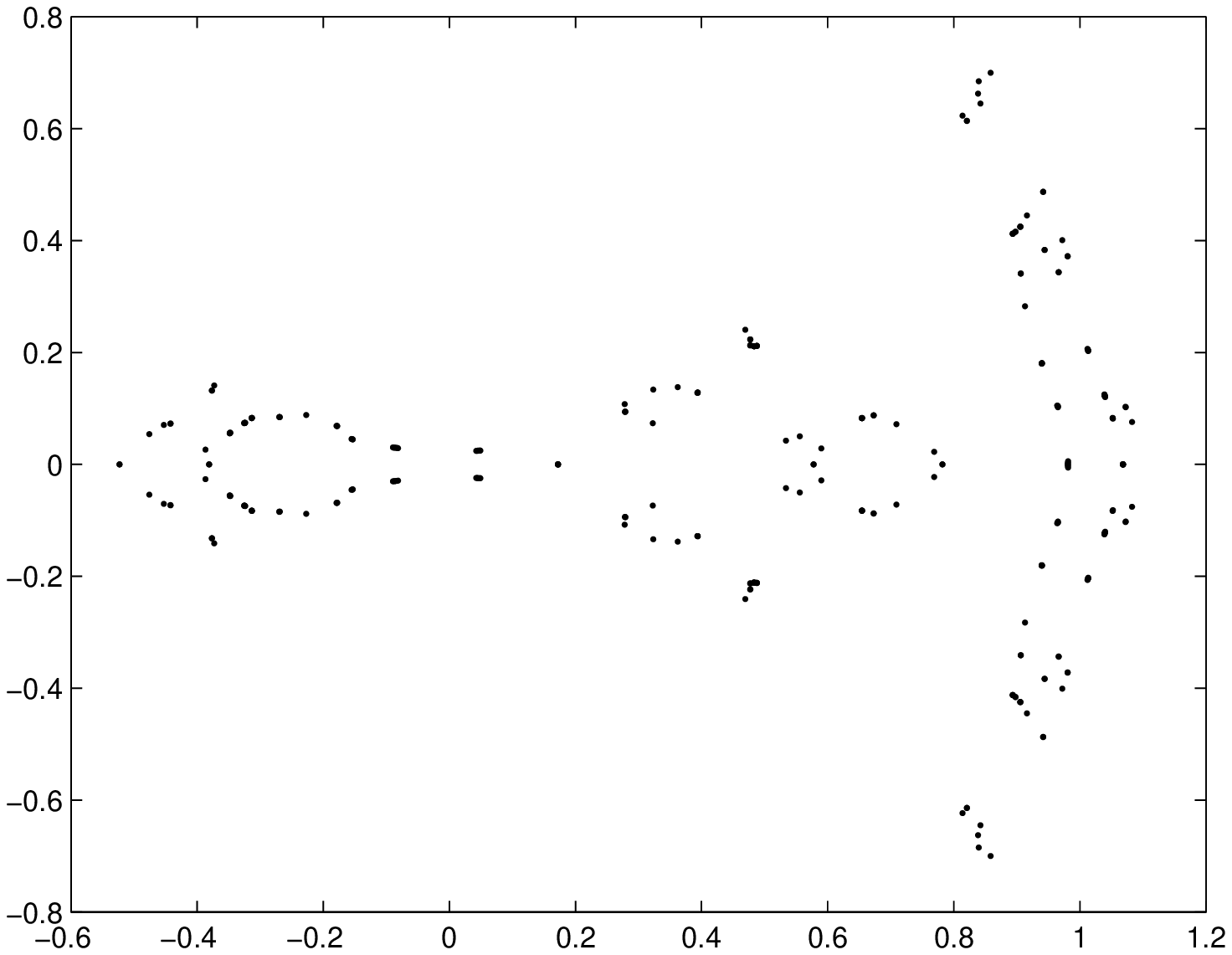}
\includegraphics[width=0.325\textwidth]{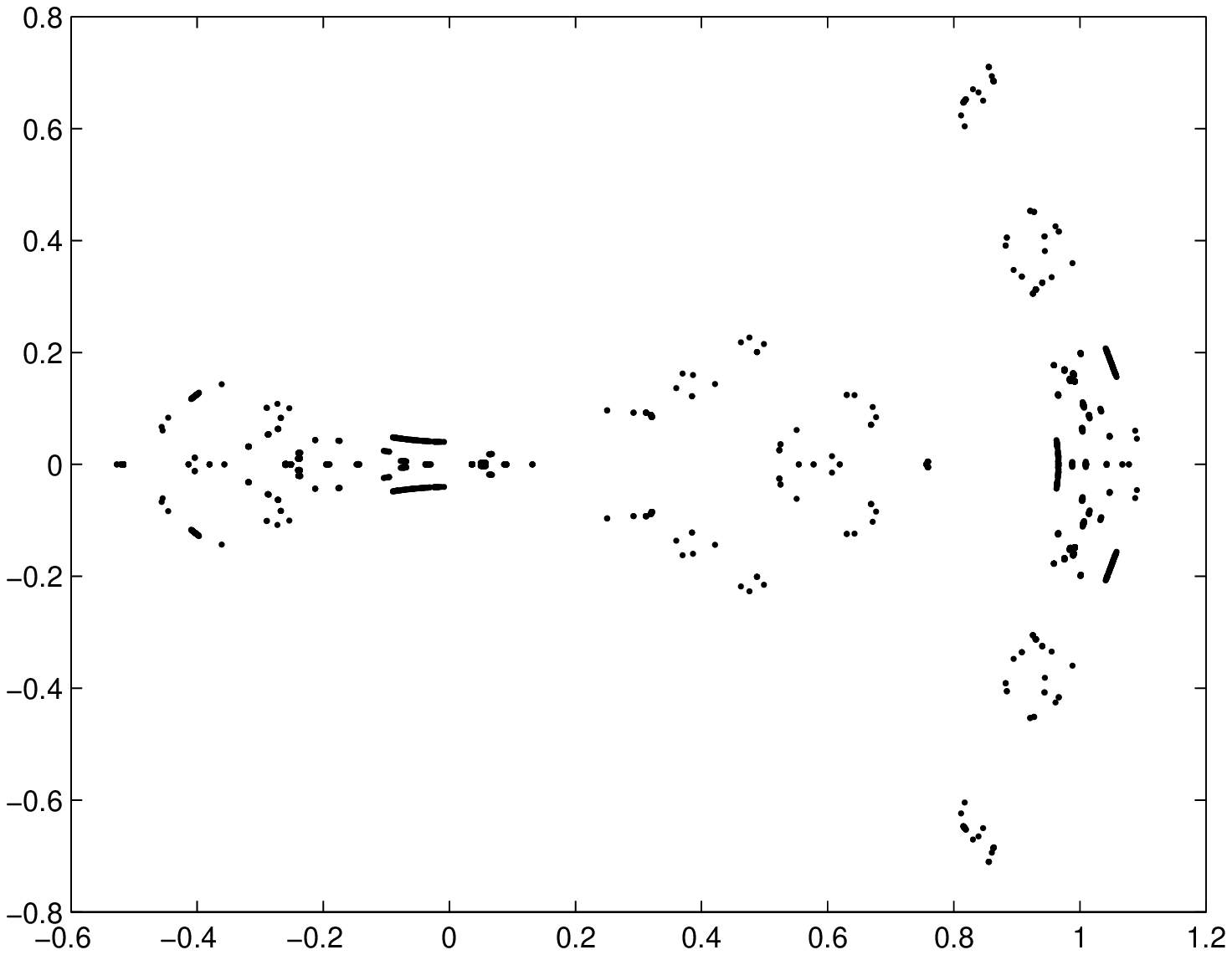}
\includegraphics[width=0.325\textwidth]{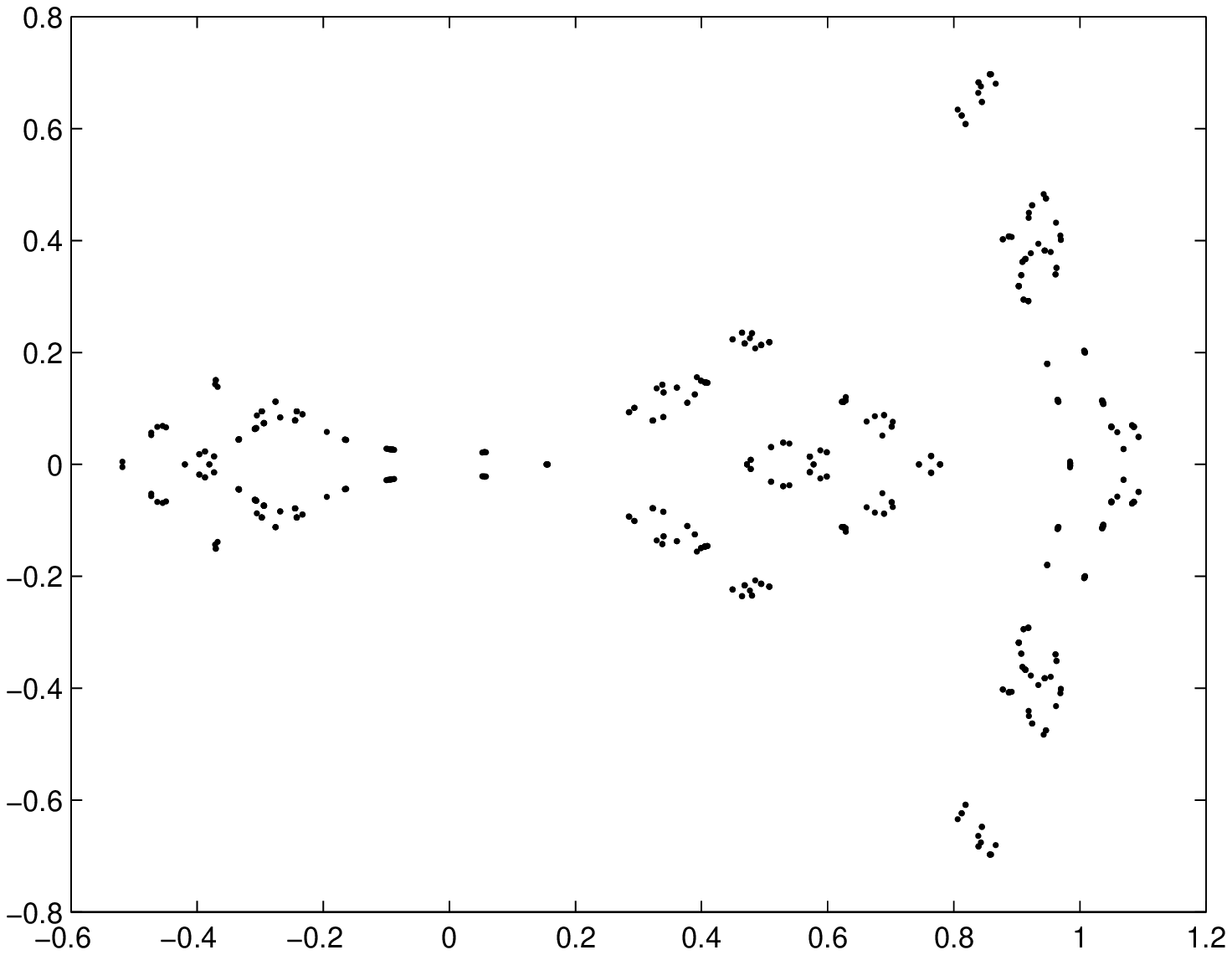}
\includegraphics[width=0.325\textwidth]{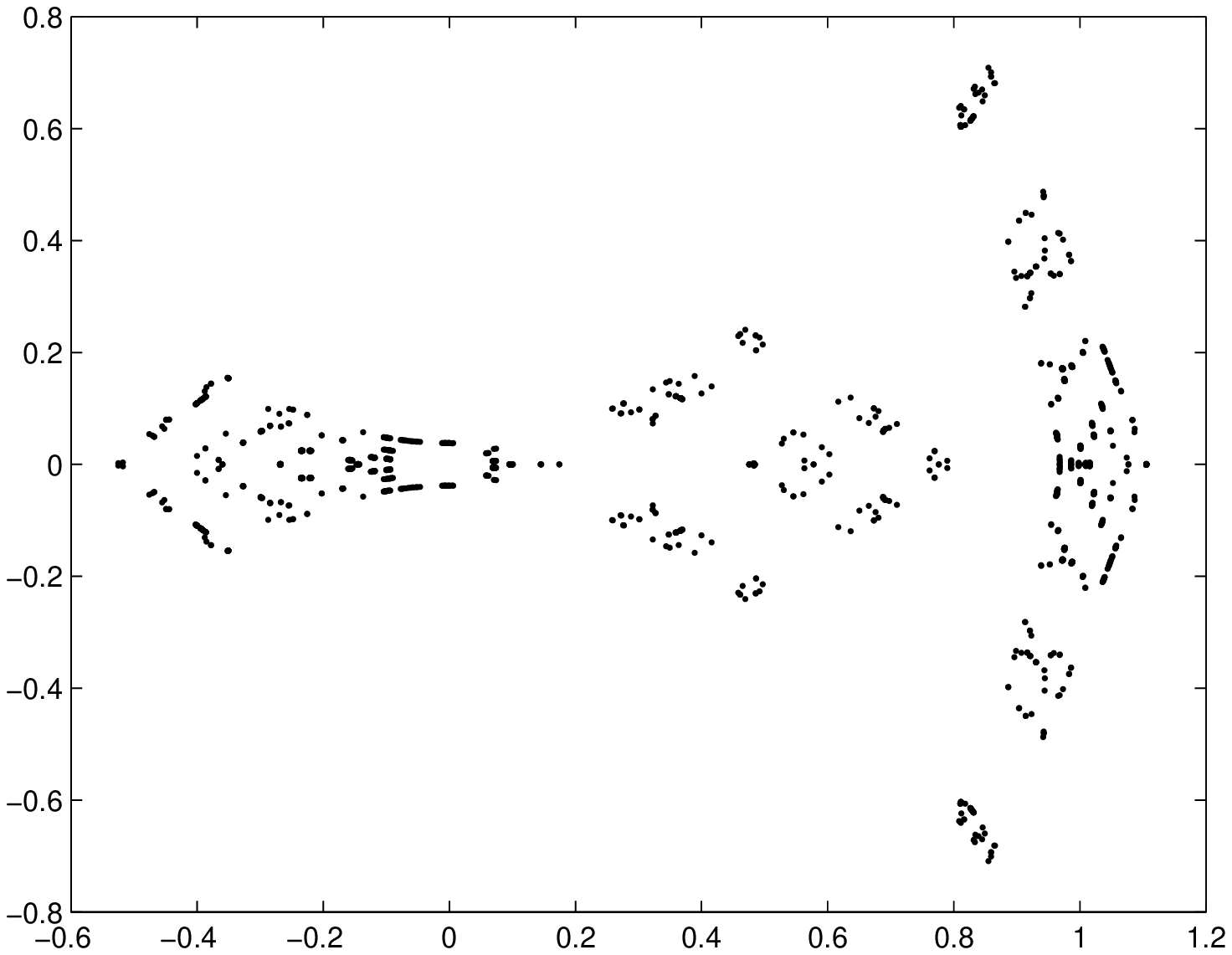}
\includegraphics[width=0.325\textwidth]{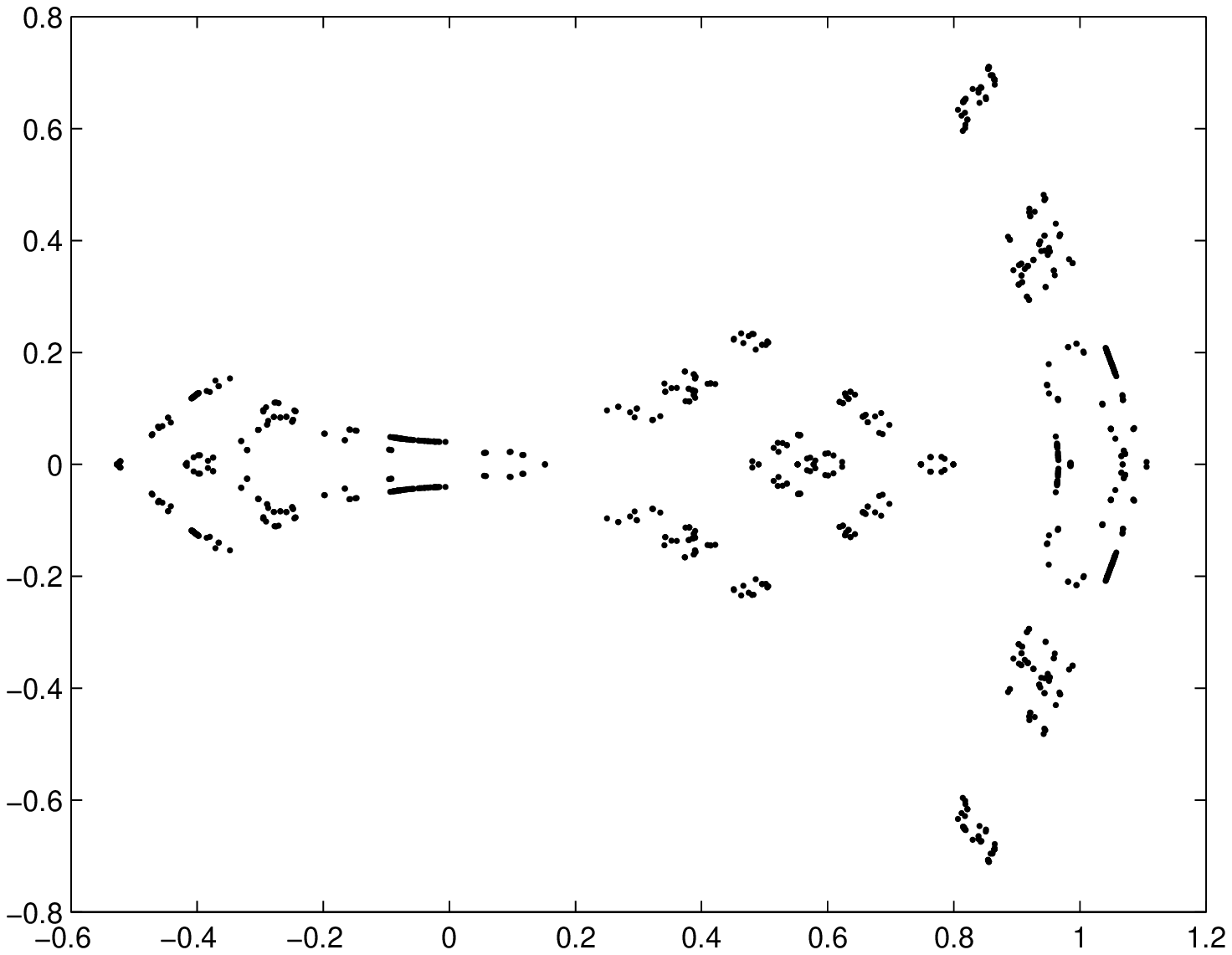}
\caption{The enclosures of the periodic cycles of $F_*$ of period up to $15$.}\label{fPerOrb}
\end{center}
\end{figure}

\begin{figure}[h]
\begin{center}
\includegraphics[width=0.48\textwidth]{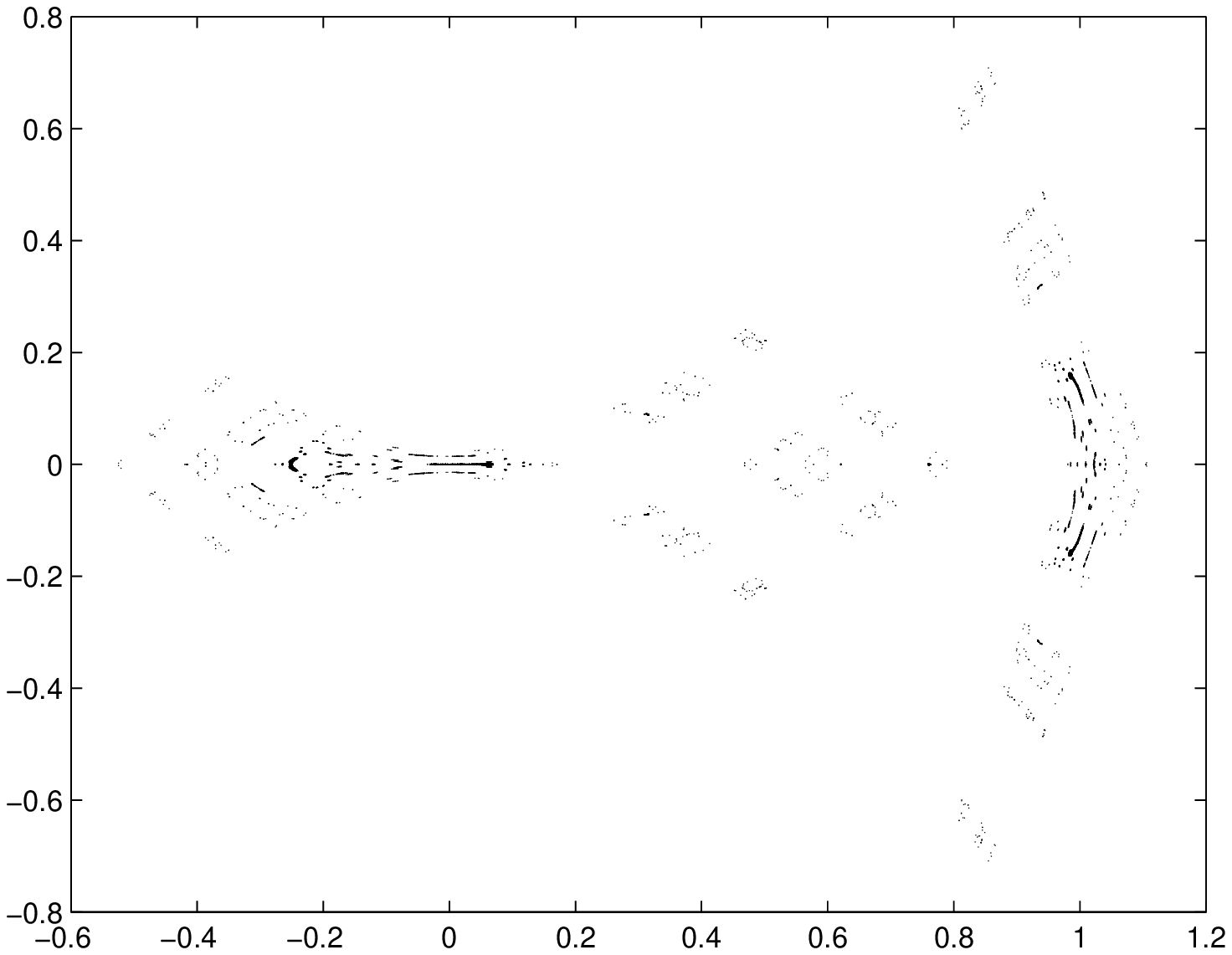}
\includegraphics[width=0.48\textwidth]{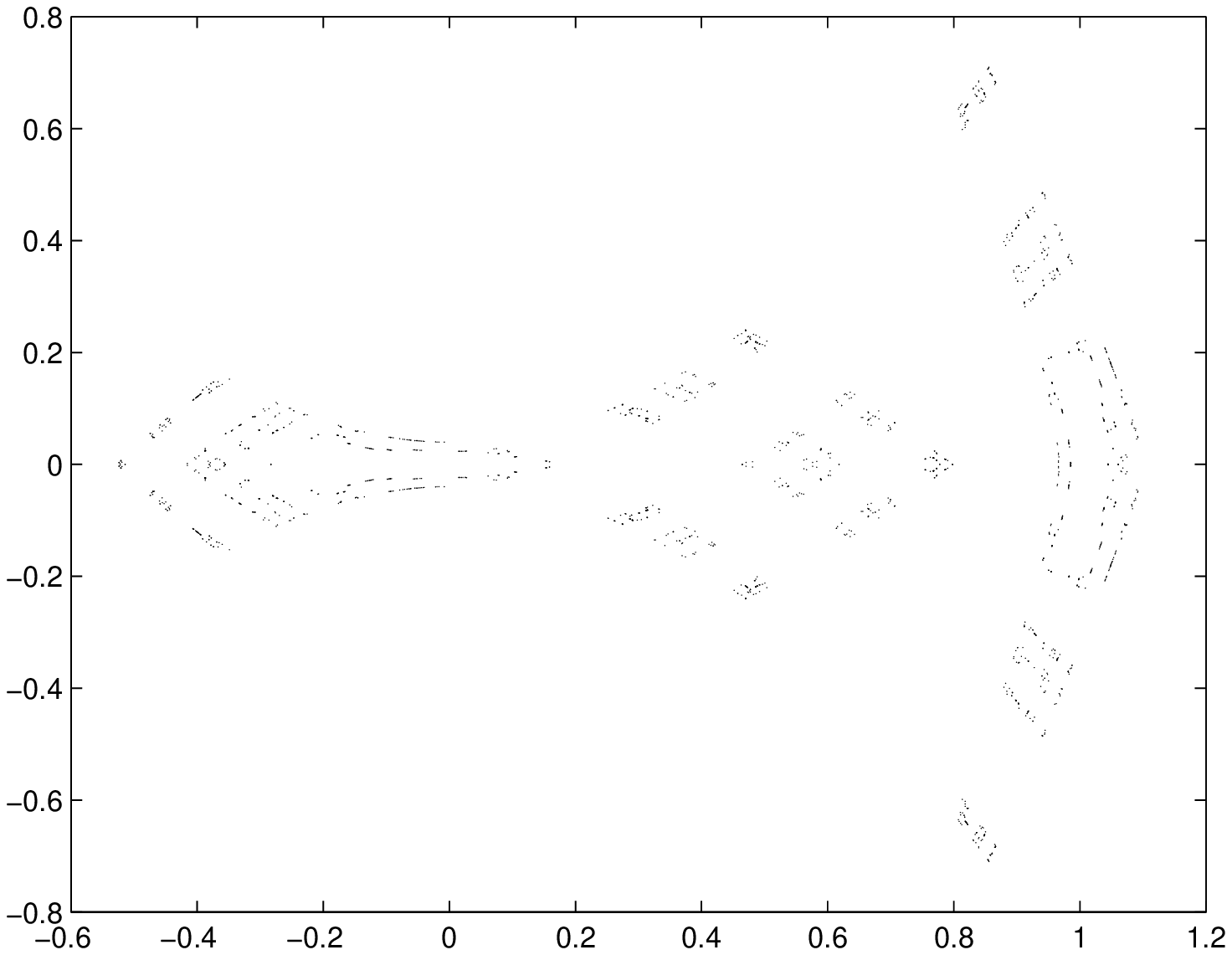}
\includegraphics[width=0.48\textwidth]{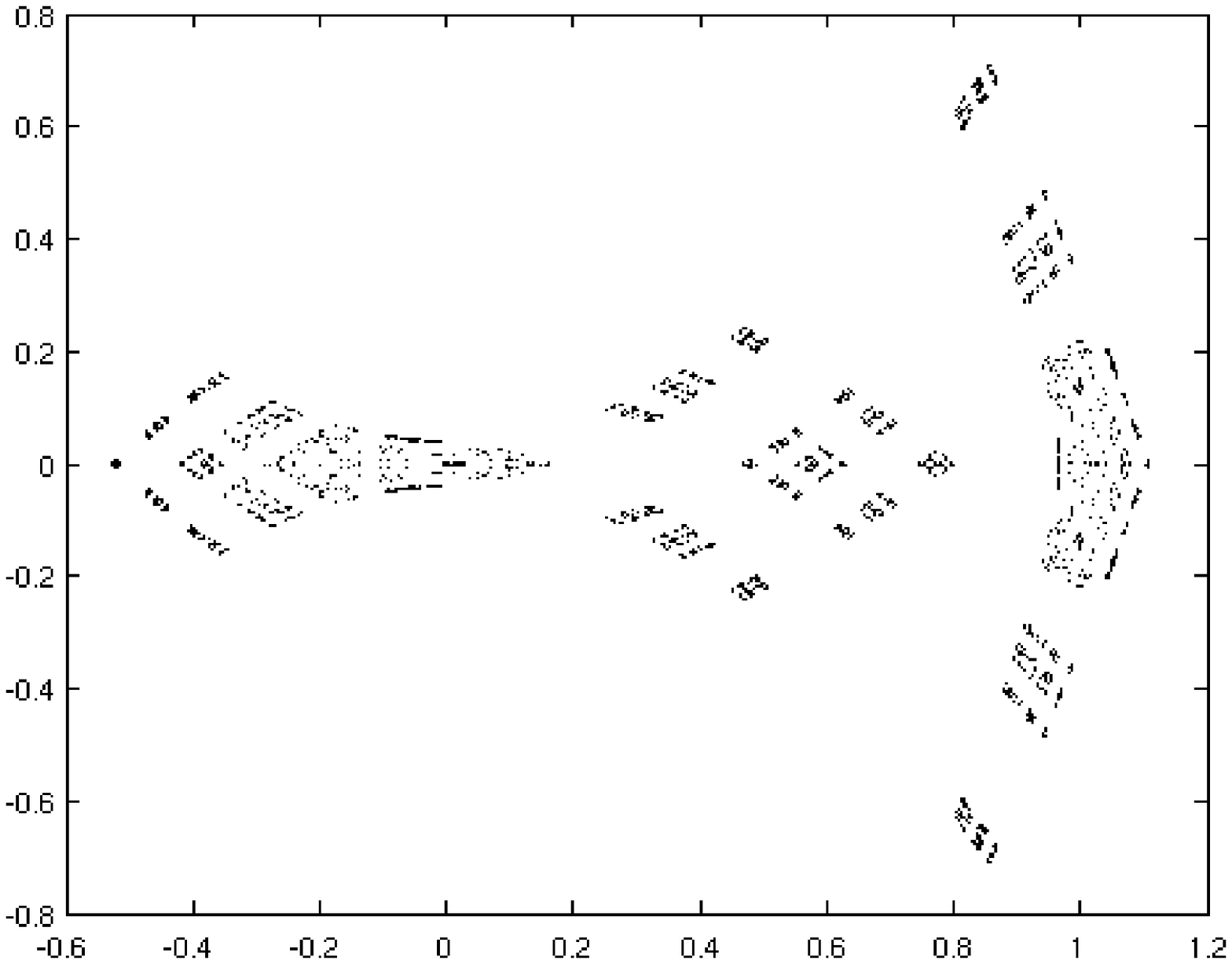}
\includegraphics[width=0.48\textwidth]{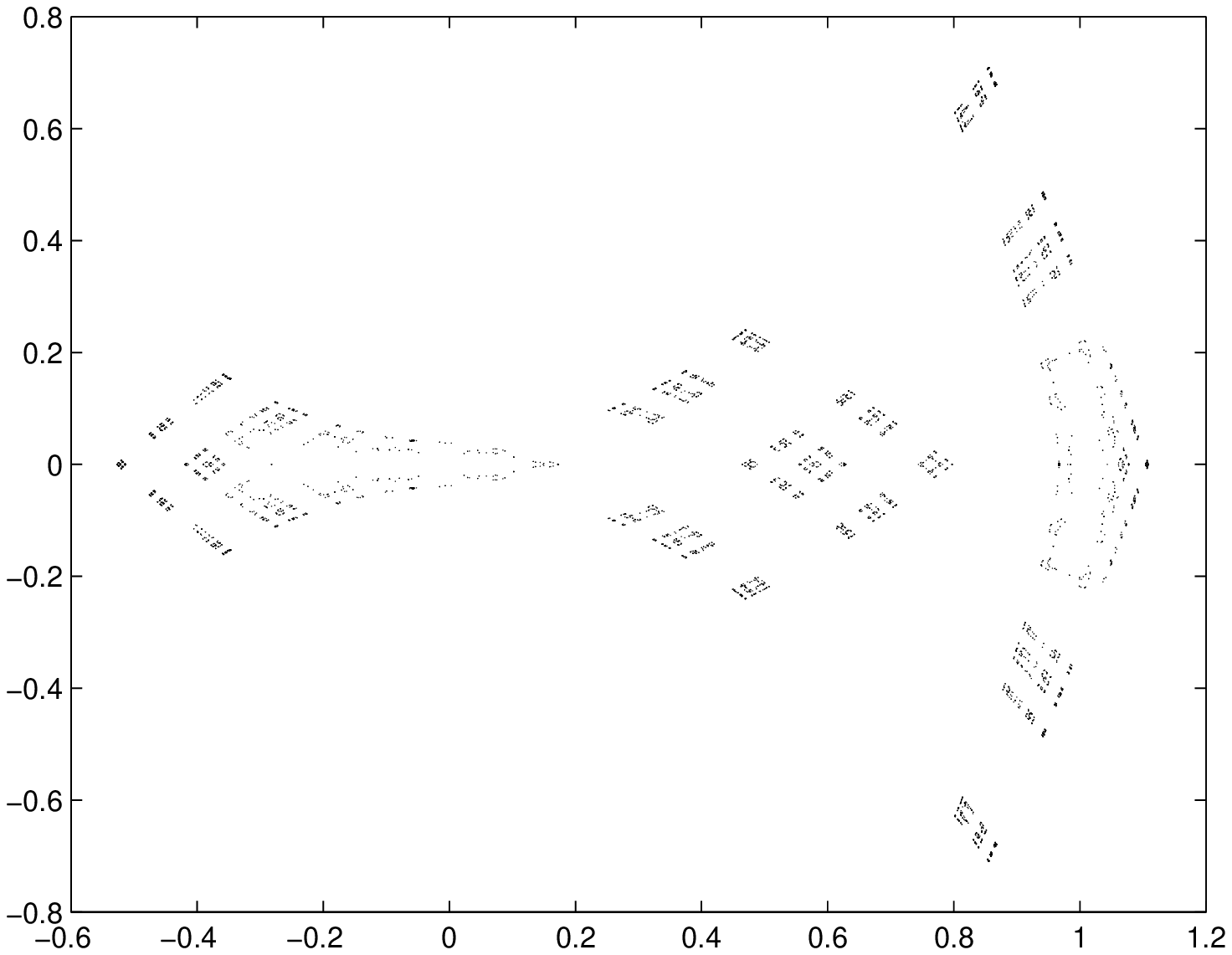}
\caption{The enclosures of the periodic cycles of $F_*$ with periods between $16$ and $19$.}\label{fPerOrbH}
\end{center}
\end{figure}

The periodic cycles with periods less than $20$ are shown in Figures \ref{fPerOrb} and \ref{fPerOrbH}. To illustrate the complexity of the computations, approximate total runtimes for the last $4$ computed periods are given in table \ref{tCPU}.

\begin{table}[h]
\begin{center}
\begin{tabular}{c|r}
$n$ & runtime  \\ \hline
$16$ & $3025$ \\
$17$ & $35$ \\
$18$ & $11060$ \\
$19$ & $1495$ \\
\end{tabular}
\caption{The total runtimes, in hours, for the proof at different periods.}
\end{center}
\end{table}\label{tCPU}

\ack
The author is funded by a postdoctoral fellowship from \textit{Vetenskapsr\aa det} (the Swedish Research Council). He would like to thank Denis Gaidashev and John Guckenheimer for useful discussions on the subject.

\end{document}